\documentclass{siamart190516}
\usepackage[utf8]{inputenc}
\usepackage{graphicx}
\graphicspath{{./figures/}}
\usepackage{amsmath}
\usepackage{amsfonts}
\usepackage{amssymb}
\usepackage{tikz}
\usepackage{colortbl}
\definecolor{lightblue}{rgb}{0.9,0.9,1}
\ifpdf
  \DeclareGraphicsExtensions{.eps,.pdf,.png,.jpg}
\else
  \DeclareGraphicsExtensions{.eps}
\fi

\newsiamremark{remark}{Remark}
\newsiamremark{hypothesis}{Hypothesis}
\crefname{hypothesis}{Hypothesis}{Hypotheses}
\newsiamthm{claim}{Claim}

\DeclareMathOperator*{\argmin}{arg\,min}
\DeclareMathOperator{\Tr}{Tr}

\headers{RBF-PUM for thin geometries and PDEs}{Larsson, Villard, Tominec, Sundin, Michael, and Cacciani}

\title{An RBF partition of unity method for geometry reconstruction and PDE solution in thin structures}

\author{Elisabeth Larsson\thanks{Uppsala  University, Department  of  Information  Technology, SE-751 05 Uppsala, Sweden
  (\email{elisabeth.larsson@it.uu.se},
    \email{ulrika.sundin@it.uu.se},
    \email{andreas.michael@it.uu.se}
  ).}
\and Pierre-Frédéric Villard\thanks{Université de Lorraine, CNRS, Inria, LORIA, Nancy, France 
(\email{pierrefrederic.villard@loria.fr}).}
\and Igor Tominec\thanks{Stockholm University, Department of Mathematics, SE-106~91 Stockholm, Sweden
(\email{igor.tominec@math.su.se}).}
\and Ulrika Sundin\footnotemark[1]
\and Andreas Michael\footnotemark[1]
\and Nicola Cacciani\thanks{Basic and Clinical Muscle Biology, Department of Physiology and Pharmacology and Center for Molecular Medicine, Karolinska Institutet, Stockholm, Sweden (\email{nicola.cacciani@ki.se}).}
\funding{This work was supported by the Swedish Research Council, grants
no.~2016-04849 and no.~2020-03488.}}

\begin{document}

\maketitle
\begin{abstract}
The main respiratory muscle, the diaphragm, is an example of a thin structure. We aim to perform detailed numerical simulations of the muscle mechanics based on individual patient data. This requires a representation of the diaphragm geometry extracted from medical image data. We design an adaptive reconstruction method based on a least-squares radial basis function partition of unity method. The method is adapted to thin structures by subdividing the structure rather than the surrounding space, and by introducing an anisotropic scaling of local subproblems. The resulting representation is an infinitely smooth level set function, which is stabilized such that there are no spurious zero level sets. We show reconstruction results for 2D cross sections of the diaphragm geometry as well as for the full 3D geometry. We also show solutions to basic PDE test problems in the reconstructed geometries.    
\end{abstract}



\section{Introduction}
This work has been performed with a particular application in mind. We aim to build a realistic model of the human respiratory system, with a special focus on the respiratory muscles, to enable in silico investigation of the adverse effects of mechanical ventilation on the muscle tissue, as well as the benefits of potential improvements to the mechanical ventilation process~\cite{invive_project}. To meet the specific challenges of this application we develop numerical tools that also have a wider use. The types of problems that we consider have the following characteristics:
\begin{itemize}
\item The problem can be formulated as (a system of) partial differential equations (PDE). Furthermore, PDE coefficents, boundary data, and, if present, inital values, are smooth and compatible such that they admit a smooth solution.
\item The PDE is formulated in a smooth, not easily parametrized geometry.
\item The geometry has a high aspect ratio. For the diaphragm, which is the main respiratory muscle, the ratio is approximately 1:100.
\end{itemize}

The simulation probem consists of two main parts. First, the geometry needs to be represented in such a way that its properties, such as boundary location and surface normal, can be incorporated into the PDE formulation. Then, the PDE is solved, to generate data of interest for the application specialist.
The assumed smoothness, present also for other biomechanical problems, suggests that we may benefit from using a high order method. Coupled with a non-trivial geometry, methods based on radial basis functions (RBF) have an advantage regarding the flexibility and ease of implementation. In this paper, we have chosen to use an unfitted least-squares RBF partition of unity method~\cite{LaShchHer17} both for the reconstruction and for the PDE solution to create a unified framework for the whole simulation problem.

In our previous work on the diaphragm application~\cite{invive,ToBre21,TVLBC21}, we  used the RBF-FD method~\cite{FoFly15book} and an unfitted least-squares RBF-FD method~\cite{ToLaHe21,ToBre21},
that has improved numerical stability properties in the presence of derivative boundary conditions compared with collocation RBF-FD. In~\cite{invive} and~\cite{TVLBC21}, linear elasticity model problems were solved and we investigated the need for resolution of the thin dimension and how to handle mixed boundary conditions in a smooth setting, respectively. In~\cite{ToBre21} the focus is on the unfitted method, Possion model problems are solved, and the diaphragm geometry is used as a demonstrator.

Reconstruction of a thin structure has many similarities with surface reconstruction, but with the added complication of interference between the two nearby surfaces~\cite{ASLGG17}. RBF-based methods (kernel methods) have been used for surface reconstruction for more than two decades~\cite{Carr01,Piret12}. As in these references, we focus on implicit surface representations, where a level set function determines which points are on the surface (zero level set), outside (positive values), or inside (negative values). Due to the typically large data sets, localized methods with low computational complexity are preferred. Popular approaches include using compactly supported RBFs~\cite{PaSka11}, using a partion of unity method combined with an adaptive oct tree subdivision of space~\cite{OhBeSe06,QiWaWu06,YWZP08}, or using quasi interpolation~\cite{LiWaBruWa16}.
Least-squares approximation, regularization, and smoothing are used to suppress noise in the data~\cite{Carr03,OhBeSe06,YWZP08,LiWaBruWa16}. Data on the surface is often complemented by exterior and interior data points to guide the direction of the level set function~\cite{Carr01,Piret12}, whereas our work is inspired by~\cite{PiDu16}, where gradient information is incorporated using a Hermite type RBF approximation. In the recent paper~\cite{DraFuWri21}, the gradient field is instead approximated by a curl-free RBF partition of unity method, which then implicitly provides the level set function. Another recent paper \cite{WMTPM21},  considers reconstruction of thin plant leaves. However, in this case the aspect ratio is high enough that the leaves are considered as surfaces. 

Our objective is to create an infinitely (or high-order) smooth geometry reconstruction in the form of an implicit surface representation, starting from a noisy point cloud representation. Furthermore, to facilitate the use of the geometry in the PDE problem, we require the level set function to be free from  spurious zero level sets inside and outside the geometry.
 We do this using the least-squares RBF partition of unity method (RBF-PUM) from~\cite{LaShchHer17}, where
node generation is decoupled from the geometry of the computational domain, with the following specific adaptions:
\begin{itemize}
\item The thin structure is adaptiviely subdivided instead of the space.
\item The same type of gradient conditions as in~\cite{PiDu16}, but without the Hermite basis, are combined  with exterior and interior stabilizing conditions.
\item The risk of interference between nearby surfaces is reduced by anisotropic scaling of the local subproblems. For a somewhat related approach, see~\cite{CaLaMoMo06}.
\end{itemize}

Least-squares RBF-PUM was applied to elliptic PDEs in~\cite{LaShchHer17}. Model problems with parametrized geometries were solved. Furthermore, a uniform grid-based distribution of spherical patches covering the domain was used. In this paper, we evaluate the performance of least-squares RBF-PUM combined with our smooth geometry reconstruction and with a non-uniform distribution of cylindrical patches.
To achieve convergence in a finite precision setting when using infinitely smooth RBFs, a stable evaluation method, such as the RBF-QR method~\cite{FoPi07,FoLaFly11,LLHF13}, is needed. 
The 3D RBF-QR implementation~\cite{LLHF13} previously did not contain all second derivative operators. For the experiments in this paper, we have updated the method with the full set~\cite{RBFQR3Dmixed}. As for the reconstruction problem, in the PDE solver, we allow anisotropic scaling of the local subproblems. However, in this case the scaling is employed to balance the resolution across the thin dimension against the resolution along the other two dimensions.

A convergence estimate for least-squares RBF-PUM was derived in~\cite{LaShchHer17}. The estimate contains a stability norm, that was investigated numerically, and the interpolation error, for which theoretical estimates are available. The convergence rate is provided by the interpolation error and we could show numerically that the stability norm for least-squares RBF-PUM stays constant under patch refinement, which is not the case for collocation RBF-PUM.
Recently, in~\cite{ToLaHe21}, we used a different technique to analyze the stability of a least-squares RBF-FD method. This allowed us to prove stability, and we plan to revisit the analysis of least-squares RBF-PUM using this technique. In this paper, we adapt the interpolation error estimates to our specific setting, to explain the observed convergence rates.

The outline of the paper is as follows. In \cref{sec:data} we describe the data sets we use for the geometry reconstruction. \Cref{sec:rbfpum} briefly introduces the least-squares RBF-PUM method. \Cref{sec:cover} and \cref{sec:weights} provide details on how we generate the partition of unity cover and weight functions. \Cref{sec:geom} describes how we use least-squares RBF-PUM for geometry reconstruction with results in \cref{sec:exp}. \Cref{sec:PDE} describes how the method is applied to the PDE problem with numerical results in \cref{sec:pdeexp}. The paper ends with conclusions in \cref{sec:conc}.

\section{The noisy input geometry data}\label{sec:data}
For our target application, the input data is extracted from 3D medical images. The extraction process is partly based on manual segmentation using physiological knowledge~\cite{villard11}, and partly automatic processing using the Marching Cube algorithm to create a mesh representing the geometry, followed by the decimation method in~\cite{levy2010} to generate the final data set. These steps are described in more detail in~\cite{invive}. The data consists of a point cloud $Y=\{\underline{y}_k\}_{k=1}^M\in\mathbb{R}^d$, for $d=2$ or $d=3$, representing the curve/surface and approximate surface normals $\underline{n}_k$, $k=1,\ldots,M$ at the corresponding points. These algorithms generate data that is quasi uniformly distributed over the surface. For simplicity, we assume this property in the algorithms that are derived in the following sections, but non-uniform distributions could also be taken into account.

In addition to the data, we make use of the initial approximate surface representation as a mesh provided by the extraction algorithm.
If an initial mesh representation was not available, we could instead create local triangulations, on demand, that would serve our purposes.  

\begin{figure}[!htb]
  \centering
  \includegraphics[width=0.405\textwidth]{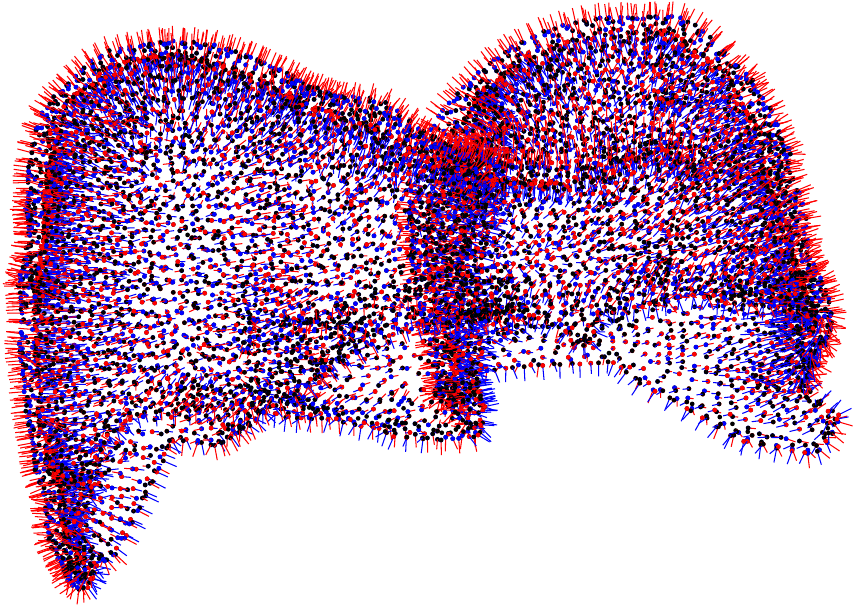}\hspace*{0.05\textwidth}
  \raisebox{0.063\textwidth}{\includegraphics[width=0.45\textwidth]{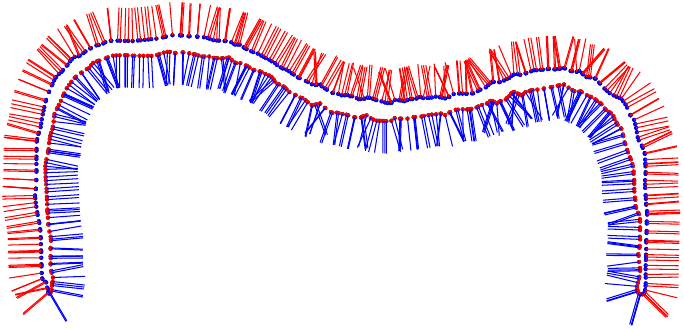}}
  \caption{The initial point cloud data in 3D (left) and 2D (right), together with the initial normal data. The surface is divided into an inner and outer part whose normals are displayed in different colors.}
  \label{fig:data}
\end{figure}

An example of a point cloud representing a diaphragm geometry is shown in the left part of~\cref{fig:data}. In the numerical experiments we also use a 2D cross section of the diaphragm. This data is shown in the right part of the figure. It is visually clear from the normals in the 2D case that the data is noisy and contains artifacts from the manual segmentation, taking the form of bumps.

\section{Least-squares RBF-PUM}\label{sec:rbfpum}
In a partition of unity (PU) method, the computational domain $\Omega\subset\mathbb{R}^d$ is covered by a set of overlapping patches $\{\Omega_j\}_{j=1}^P$, such that $\bar{\Omega}\subset\cup_{j=1}^P\bar{\Omega}_j$. We require each point $\underline{x}=(x_1,\ldots,x_d)\in\bar{\Omega}$ to be in the interior of at least one patch. We also require the overlap to be mild in the sense that no more than $K$ patches overlap at any given point. We construct a set of non-negative partition of unity weight functions $\{w_j(\underline{x})\}_{j=1}^P$, where $w_j$ is compactly supported on $\Omega_j$. The partition of unity property implies that
$\sum_{j=1}^P w_j(\underline{x}) = 1$, for $\underline{x}\in\cup_{j=1}^P\Omega_j$.

A global PU approximation $\tilde{u}(\underline{x})$ is constructed from local approximations $\tilde{u}_j(\underline{x})$ on each patch and the weight functions as
\begin{equation}
  \tilde{u}(\underline{x}) = \sum_{j=1}^Pw_j(\underline{x})\tilde{u}_j(\underline{x}).
  \label{eq:rbfpum}
\end{equation}  
The smoothness of the global approximation is determined by the smoothness of the weight functions and of the local approximations. In this paper, we use RBF-PUM both for the geometry representation and for solving a PDE problem formulated in that geometry. For the PDE problem it is enough if the smoothness of the weight function supports the differential operators of the PDE. However, to construct a high-order method for the PDE problem, we need high-order smoothness of the geometry. That is, for the geometry representation, we require high-order smoothnes also of the weight functions.

For the local RBF approximations, we make use of a template node layout in a reference patch $\Omega_0$, see~\cite{LaShchHer17}. The patches $\Omega_j$ are allowed to have different sizes, different aspect ratios, and different orientations. For each patch, we introduce a map $T_j\,:\,\Omega_j\mapsto \Omega_0$ that takes points in $\Omega_j$ into the reference patch. We let the local approximations be given by
\begin{equation}
  \tilde{u}_j(\underline{x}) = \tilde{v}_j(T_j(\underline{x})) \triangleq \tilde{v}_j(\underline{x}^\prime),
  \label{eq:uofv}
\end{equation}
where $\tilde{v}_j(\underline{x}^\prime)$ is an RBF approximation in the reference patch, of the form
\begin{equation}
  \tilde{v}_j(\underline{x}^\prime) = \sum_{i=1}^{n}\lambda^j_i\phi(\|\underline{x}^\prime-\underline{x}^\prime_i\|)\triangleq \sum_{i=1}^{n}\lambda^j_i\phi_i(\underline{x}^\prime),
  \label{eq:RBF}
\end{equation}
where $\lambda_i^j$ are unknown coefficients, $\phi(r)$ is an RBF, and $X^\prime=\{\underline{x}^\prime_i\}_{i=1}^n$ is the set of RBF center points in the reference patch, with fill distance $h^\prime$.
The corresponding points in $\Omega_j$ are $\underline{x}_i^j=T_j^{-1}(\underline{x}_i^\prime)$, $i=1,\ldots,n$, with fill distance $h$. We let $\underline{u}_j=(\tilde{u}_j(\underline{x}_1^j),\ldots,\tilde{u}_j(\underline{x}_n^j))^\intercal =(\tilde{v}_j(\underline{x}_1^\prime),\ldots,\tilde{v}_j(\underline{x}_n^\prime))^\intercal$ and $\underline{\lambda}_j=(\lambda_1^j,\ldots,\lambda_n^j)^\intercal$. Then, using~\cref{eq:RBF}, we have $\underline{u}_j=A\underline{\lambda}_j$, where the local interpolation matrix $A$ has elements $a_{ki}=\phi(\|\underline{x}_k^\prime-\underline{x}^\prime_i\|)$. For commonly used infinitely smooth RBFs such as the multiquadric with $\phi(r)=\sqrt{1+r^2}$, the Gaussian with $\phi(r)=e^{-r^2}$, and the inverse multiquadric with $\phi(r)=1/\sqrt{1+r^2}$, the interpolation matrix $A$ is non-singular~\cite{Scho38,Micchelli86}. This implies that we can safely use $\underline{\lambda}_j=A^{-1}\underline{u}_j$ in~\cref{eq:RBF} to get the local RBF approximations expressed in terms of the nodal values as
\begin{equation}
  \tilde{v}_j(\underline{x}^\prime) = (\phi_1(\underline{x}^\prime),\ldots,\phi_n(\underline{x}^\prime))A^{-1}\underline{u}_j.
  \label{eq:RBF2}
\end{equation}
Note that in this formulation of RBF-PUM, each local approximation has its own nodal values, also in the overlapping regions. The total number of nodal unknowns to determine is $N=nP$.

When we use RBF-PUM to approximate a function, the local subproblems can be solved independently, and are then blended into the global approximation. For a PDE, we instead need to assemble a global system of equations to get a well-posed discretized problem. We describe a generic time-independent linear PDE by
\begin{equation}
  \mathcal{L}(\underline{x})u(\underline{x}) = f(\underline{x}),\quad \underline{x}\in\bar{\Omega}.
\end{equation}
where $\mathcal{L}(x)$ represents the differential operators of the PDE and its boundary conditions, and $f(x)$ represents all right hand side functions. To form the global least-squares system of equations, we sample the PDE problem at the points  $\{\underline{y}_k\}_{k=1}^M$, where $M>N$. If we let $\mathcal{L}_k=\mathcal{L}(\underline{y}_k)$, then the equations are given by
%
\begin{equation}
  \mathcal{L}_k\tilde{u}(\underline{y}_k) = \sum_{j=1}^P\mathcal{L}_k(w_j(\underline{y}_k)\tilde{u}_j(\underline{y_k}))
  =f(\underline{y}_k).
  \label{eq:PDERBFPUM}
\end{equation}
Practically, the matrix is assembled by adding the contributions from each patch. We write the global system as $L\underline{u}=\underline{f}$, where the vector $\underline{u}=(\underline{u}_1^\intercal,\ldots,\underline{u}_P^\intercal)^\intercal$ and the vector $\underline{f} = (f(\underline{y}_1),\ldots,f(\underline{y}_M))^\intercal$. Note that we differentiate the product of the weight functions and the local approximations such that each differential operator is expanded into several terms.

When solving the least squares system of equations, the scaling of the equations corresponding to different operators becomes important for convergence. Here, we use the insights from~\cite{ToLaHe21}, where we dereived and analyzed a least-squares RBF-FD method. First we group the equations according to the type of operator. If there are $M_i$ equations of type $i$, then we scale these equations with $1/\sqrt{M_i}$. With this scaling, the discrete least-squares problem is a discrete quadrature approximation of the corresponding continuous least-squares problem. Then, we increase the weight of the lower order boundary equations by multiplying these with a small negative power of the fill distance $h$ of the reference node set $X$. This improves stability without affecting the rate of convergence.

\section{Adaptive cover generation for thin structures}\label{sec:cover}
\begin{figure}[!htb]
  \centering
  \includegraphics[width = 0.45\textwidth]{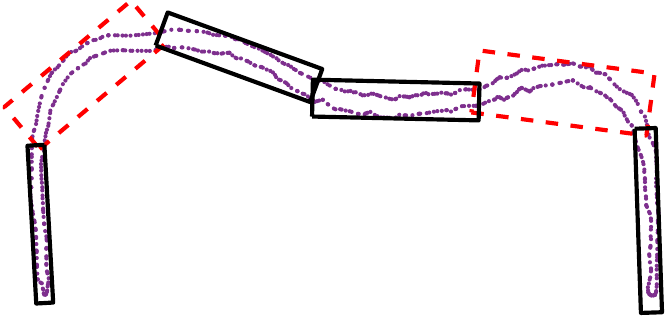}
  \includegraphics[width = 0.45\textwidth]{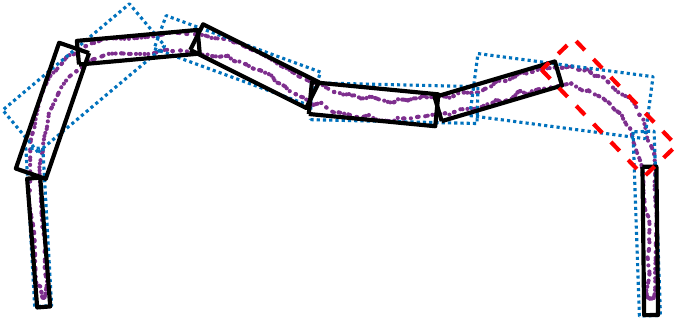}
  \includegraphics[width = 0.45\textwidth]{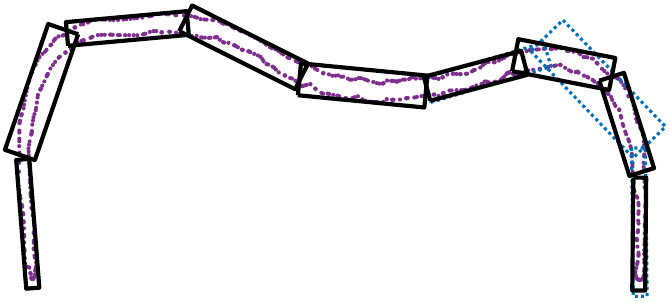}
  \includegraphics[width = 0.441\textwidth]{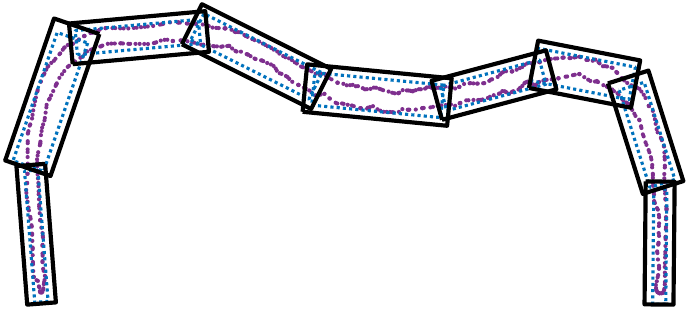} 
  \caption{
    The adaptive cover generation with the initial patches (top left), the intermediate result (top right), and the final non-overlapping patches (bottom left). In each step, patches that that are marked for refinement (dashed lines) and the previous iteration (dotted lines) are indicated. The data points (dot markers) are shown for reference. Finally, the resulting patches after ensuring an overlap $\delta\approx0.15$ are shown (bottom right). }
\label{fig:patch2D}  
\end{figure}

\begin{figure}[!htb]
\centering
  \includegraphics[width = 0.45\textwidth]{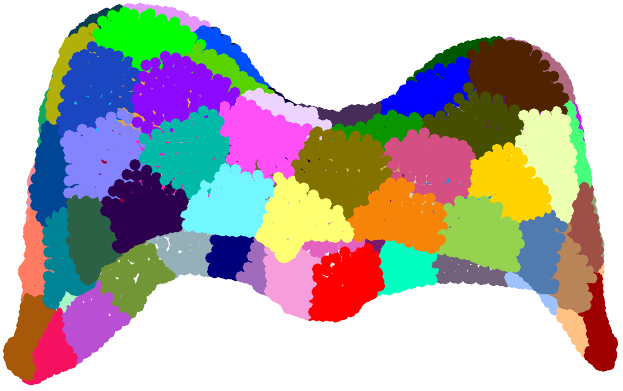}
  \includegraphics[width = 0.42\textwidth]{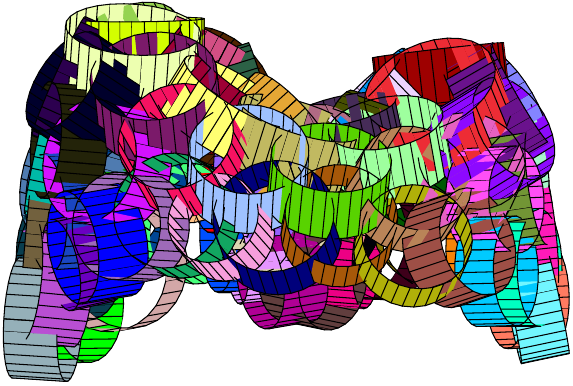}
  \caption{The result of $k$-means clustering of 10\,276 nodes into 75 clusters (left) and the $P=79$ final patches generated by the 3 steps of the cover algorithm (right). Only the mantle of the cylindrical patches is drawn to make the intersections visible.}
  \label{fig:patch3D}
\end{figure}

%
Common choices of patch shapes for RBF-PUM in the literature are discs and spheres or their anisotropic counterparts, see, e.g.,~\cite{SVHL15,LaShchHer17,CheShch18,Cav19,Cav20}. These work well when the computational domain is isotropic or has a fixed direction of anisotropy. Here, we are considering thin curved structures, and it becomes relevant to use a patch shape that adheres to this property. We have chosen to use cylindrical patches that are radial in directions tangential to the surface and with a height that is related to the thickness of the domain. We also make the assumption that there is only one layer of patches in the thickness direction. We let the reference patch, centered at the origin, have radius $R_0$ and height $H_0$. An example of two-dimensional cylindrical patches can be seen in~\cref{fig:patch2D} and three-dimensional cylindrical patches are shown in~\cref{fig:patch3D}.

The starting point for the adaptive cover generation is the input data $Y$ with $M$ points and the user defined parameter $P_0$ corresponding to the initial number of patches. After the adaptive refinement, the final number of patches is $P\geq P_0$. The cover algorithm consists of the following steps: 
\begin{enumerate}
\item Generate and adaptively refine patches until a volume criterion is met.
  \begin{enumerate}  
  \item Use $k$-means clustering to form non-overlapping subsets $Y_j=\{\underline{y}_i^j\}_{i=1}^{m_j}$ of the input data with corresponding patch centers $\underline{C}_j$. 
  
  \item Use principal component analysis (PCA) on each $Y_j$ to find the patch radius $R_j$, the vertical range $[Z_{j,1}\, Z_{j,2}]$, and the transformation $T_j$ to the reference patch.
  \item Evaluate the volume criterion, flag patches for refinement that fail the criteria, and return to 1(a) if not done. 
  \end{enumerate}
\item Ensure that each boundary segment or surface triangle is entirely covered.
\item Ensure an overlap $\delta$ between the patches.
\end{enumerate}  
The volume criterion that we use in Step 1 aims to have a similar resolution (fill distance) in each patch. This also means that we refine where the curvature of the domain is large within a patch, since this increases the volume. We denote the initial mean volume by $\bar{V}_0$ and the initial sample standard deviation by $s_{V_0}$, and set the (absolute) target volume $V^*$ as
\begin{equation}
  V^* = \bar{V}_0 + \max(\bar{V}_0/10,s_{V_0}),
  \label{eq:crit1}
\end{equation}
where the constant 10 is chosen to prevent excessive refinement when patches are already similar in volume. 
When patches for which $V_i>V^*$ have been flagged for refinement, we first add also the neighbours of each of these patches to the flagged set. Then for each connected set of $p$ flagged patches, we take the union of the corresponding $p$ node subsets, change the required number of patches to $p+1$ for this subset of the data and construct these $p+1$ new patches locally. This process is illustrated by the first tree panels of~\cref{fig:patch2D}.

The first step of the patch construction for a node set is $k$-means clustering. The algorithm aims to create clusters that are of almost equal size and with a minimal sum of square distances to the centroid. An example of a clustering is shown in the left part of~\cref{fig:patch3D}. The initial patch centers $\underline{C}_j$ are taken as the centroids of the clusters.

In the local patch coordinate system, we let $\underline{C}_j$ be the origin. We let $\tilde{Y}_j=\{\underline{y}_i^j-\underline{C}_j\}_{i=1}^{m_j}$, and we denote the matrix of size $m_j\times d$ containing the points by $[\tilde{Y}_j]$. A principal component analysis of the data amounts to finding the eigenvalues $\mu_1\geq \cdots \geq \mu_d$ of $[\tilde{Y}_j]^\intercal [\tilde{Y}_j]$ and the corresponding eigenvectors $q_1,\ldots,q_d$. The first $d-1$ principal components (eigenvectors) represent the main components of the data, i.e., the tangent (plane) direction(s), while the $d$th eigenvector is the direction normal to the tangent plane.

We form the orthogonal matrix $Q_j$ with the eigenvectors as columns.
We transform the local points into the patch coordinate system and get $\hat{Y}_j=\{\underline{\tilde{y}}_i^jQ_j\}_{i=1}^{m_j}$, which allows us to compute the radius and vertical range of the patch:
\begin{align}
R_j &= \max_i\textstyle\left(\sum_{k=1}^{d-1}(\underline{\hat{y}}^j_{i,k})^2\right)^\frac{1}{2},\\
Z_{j,1} &= \min_i \underline{\hat{y}}^j_{i,d},\\
Z_{j,2} &= \max_i\underline{\hat{y}}^j_{i,d}.
\end{align}
We denote the center of the patch in the local coordinate system by $\underline{Z}_j$, which has $d-1$ zero components and $\frac{1}{2}(Z_{j,1}+Z_{j,2})$ in the last position. We also define the diagonal scaling matrix $S_j$ with diagonal elements $s_{ii}=\frac{R_0}{R_j}$ for $i=1,\ldots,d-1$ and $s_{dd}=\frac{H_0}{H_j}$. The local points in the reference patch coordinate system then become $Y_j^\prime=\{(\underline{\hat{y}}_i^j-\underline{\bar{Z}}_{j})S_j\}_{i=1}^{m_j}$. We now have all the steps to write down the transformation from the global coordinate system to the reference patch for the points in patch $j$
\begin{equation}
  T_j(\underline{x}) = ((\underline{x}-\underline{C}_j)Q_j-\underline{Z}_j)S_j=\underline{x}^\prime.
  \label{eq:Tj}
\end{equation}
When expressing conditions on functions in the physical domain in terms of the approximations in the reference patch, we need to differentiate this map. If we define the gradient as a row vector, then for a function $f(\underline{x})=g(T_j(\underline{x}))=g(\underline{x}^\prime)$ 
\begin{eqnarray}
    \nabla f & =& \nabla^\prime g S^\intercal Q^\intercal,\nonumber\\[5pt]
  \mathcal{H}_f &=& QS\mathcal{H}^\prime_g S^\intercal Q^\intercal ,  \label{eq:transder}
\\[-3pt]
  \Delta f & =& \Tr(S\mathcal{H}_g^\prime S^\intercal) 
  =\sum_{i=1}^ds_{ii}^2 \frac{\partial^2g}{\partial x_i^{\prime 2}},\nonumber
\end{eqnarray}  
where $\mathcal{H}$ denotes the Hessian matrix, primed operators refer to the reference coordinate system, and the Laplacian has been simplified using rotational invariance.

The patches that have been generated so far cover all data points, but to ensure a cover of the whole geometry, we also need to cover the potential gaps between the data points. Here we use the initial approximate surface representation, see~\cref{sec:data}. We extend the patches such that each surface element is covered. The surface elements are planar objects, while patches are convex regions. This implies that if all corners of an element are covered by a patch, then the whole element is covered. For elements with their corners in different patches, we divide the element into $d$ Voronoi regions associated with the $d$ corners, and then extend the closest patch that covers a corner to cover also the Voronoi region of that corner. When all boundary elements are covered, we have ensured that the whole initial object is covered. 

The final step of the cover algorithm ensures that there is sufficient overlap between the patches. The overlap is important because it affects the steepness of the weight functions. Across an overlap region, the involved weight functions need to transition smoothly from 0 to 1.
To handle this practically, we require that all data points have at least a at relative distance $\delta_0$ from the inside of the patch to any patch boundary.
For each point where this is not satisfied, we extend the patch that requires the smallest volume increase for including that point. The resulting relative overlap will not be precisely $\delta_0$, but close enough as long as the data points are dense enough.
The result of the overlap algorithm is shown in the last panel of~\cref{fig:patch2D} for the two-dimensional diaphragm geometry.

\begin{remark}
  If the surface data is too sparse in relation to the thickness of the geometry, the clustering can be performed only on the inner or outer surface in order to acheive the desired result with one layer of patches. We have implemented an automatic approach to label the edge, the inner surface, and the outer surface.
\end{remark}

\section{Tensor product partition of unity weight functions}\label{sec:weights}
The partition of unity weight functions are tightly connected with the cover. To build the weight functions, we start from a set of non-negative generating functions supported on the respective patches. The cylindrical patch shape that we have chosen suggests a tensor product form. Let $\psi_0(r)$ be a positive radial function supported on $[0,\,1)$, and let $\|\underline{x}\|_r =\sqrt{x_1^2+\cdots+x_{d-1}^2}$. Then the product
\begin{equation}
  \psi(\underline{x}^\prime)=\psi_0\left(\frac{\|\underline{x}^\prime\|_r}{R_0}\right)\psi_0\left(\frac{2|x_d|}{H_0}\right)
\end{equation}
is compactly supported on the reference patch with radius $R_0$ and height $H_0$. The generating function for patch $j$ is then given by
\begin{equation}
  \psi_j(\underline{x}) = \psi(T_j(\underline{x})) = \psi(\underline{x}^\prime). 
\end{equation}
The generating functions we use in this paper are the infinitely smooth bump function
\begin{equation}
  \psi_0(r) = \exp\left(-\frac{1}{(1-r^2)_+}\right),
  \label{eq:bump}
\end{equation}  
for the geometry reconstruction, and the $C^2$ Wendland function~\cite{Wend95}
\begin{equation}
  \psi_0(r) = (4r+1)(1-r)^4_+,
  \label{eq:wend}
\end{equation}  
for the PDE solutions.

\begin{figure}[!htb]
  \centering
  \begin{tikzpicture}
    \begin{scope}[shift={(0,0)}]
      \node at (-1,0.3) {\includegraphics[height=0.2\textwidth]{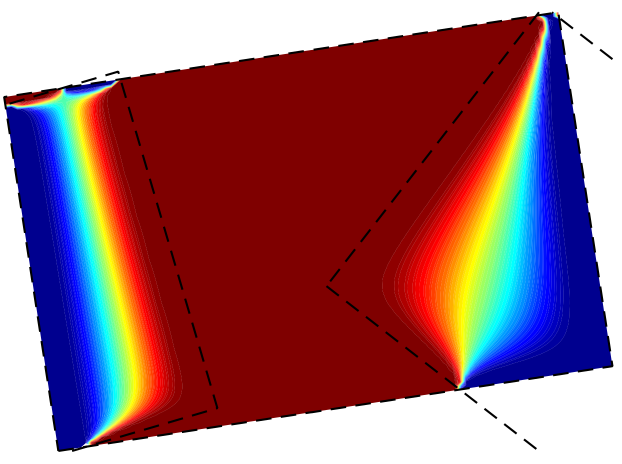}};
      \node at (3.7,0) {\includegraphics[height=0.3\textwidth]{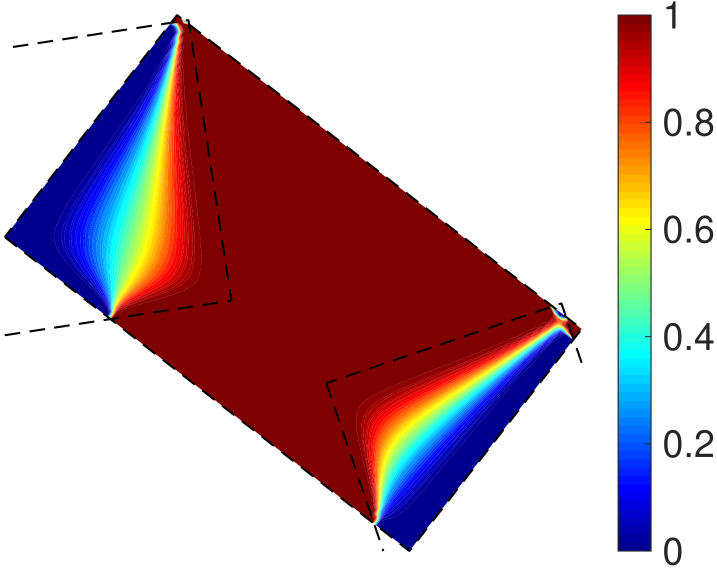}};
    \end{scope}
  \end{tikzpicture}
  \caption{Contour plots of the weight functions of two adjacent 2D patches.}
  \label{fig:weight}
\end{figure}

To construct the weight functions from the generating functions, we use Shepard's method~\cite{Shepard68}. We first define the sum
\begin{equation}
  s(\underline{x})=\sum_{j=1}^P\psi_j(\underline{x}),
\end{equation}
which is positive for all $\underline{x}\in \cup_{j=1}^P\Omega_j$. Even though the sum is over all patches, it never has more than $K$ non-zero terms at any given point due to the overlap condition. The weight function $w_j$ supported on the patch $\Omega_j$ is given by
\begin{equation}
w_j(\underline{x}) = \frac{\psi_j(\underline{x})}{s(\underline{x})}.
\end{equation}
Using a recursive differentiation rule for rational functions, we furthermore have
\begin{align}
\nabla w_j &= \frac{\nabla \psi_j -w_j\nabla s}{s},\nonumber\\
\Delta  w_j &= \frac{\Delta \psi_j-2\nabla w_j\cdot\nabla s-w\Delta s}{s},\label{eq:weightder}\\
\mathcal{H}_{w_j} & = \frac{\mathcal{H}_{\psi_j}-\nabla w_j^\intercal\nabla s
-\nabla s^\intercal\nabla w_j -w\mathcal{H}_s}{s}.\nonumber
%
%
\end{align}
An example with two weight functions is shown in~\cref{fig:weight}. Within the regions where there is no overlap, the local weight function is identically 1 and all derivatives are zero. The formulas above are therfore only needed in the overlap regions. At the outer boundary of $\cup_{i=1}^P\Omega_j$, at the points where two (or more) patches intersect, the weight functions are discontinuous, since the effective overlap at these points is 0. In the interior the weight functions have the full smoothness of the generating function.

\section{The smooth implicit geometry representation}\label{sec:geom}
We represent the geometry of the thin structures we are studying implicitly, through a level set function $\ell(\underline{x})$ that is zero at the object surface. The gradient of the level set function at the object surface provides the normal direction, but we only require $\ell(\underline{x})$ to approximately represent the distance to the object surface with negative values inside and positive outside. 

To construct $\ell(\underline{x})$ we use an RBF-PUM approximation of form~\eqref{eq:rbfpum}
\begin{equation}
  \ell(\underline{x}) = \sum_{j=1}^Pw_j(\underline{x})\ell_j(\underline{x}) =
  \sum_{j=1}^Pw_j(\underline{x})\gamma_j(T_j(\underline{x})),
\end{equation}  
where $\gamma_j$ are functions of form~\eqref{eq:RBF}--\eqref{eq:RBF2} defined over the reference patch. The least squares approximation problems for the functions $\gamma_j(\underline{x})$ are solved independently. By using the infinitely smooth multiquadric RBF for the local approximation and the infinitely smooth weight function~\cref{eq:bump}, the geometry reconstruction also becomes infinitely smooth. We need to choose the number of basis functions $n$ in the reference patch such that we ensure oversampling in all patches. The $m_j$ data points in patch $j$ generate $2m_j$ conditions based on the locations and normal directions. This provides the bound
\begin{equation}
  n \leq \min_j 2m_j. 
\end{equation}
We are in fact adding some extra conditions as described below, but we keep this bound since it is the limit when considering only the data.

The local subproblems contain three types of conditions. The conditions are formulated in the physical domain, but implemented in the reference patch, where all local problems are solved. The first condition is the zero value at the object surface
\begin{equation}
  \ell_j(\underline{y}_i^j) = \gamma_j(T_j(\underline{y}_i^j)) = 0, \quad \underline{y}_i^j\in Y_j^\rho,
  \label{eq:geomcond1}
\end{equation}
where $Y_j^\rho$ is an extension of the local data set $Y_j$ with the data points of distance at most $\rho h$ outside the patch. By adding the extra conditions outside we improve the approximation at the intersection of the radial patch boundary and the object.

The second set of conditions guides the gradient of the level set function. We follow the approach in~\cite{PiDu16} and require the gradient to have a unit component in the initial normal direction. However, since we are working with local patch approximations and not stencils, we enforce these conditions at all data points instead of only at the stencil center as was done in~\cite{PiDu16}. We use~\eqref{eq:transder} to transform the gradient to get
\begin{equation}
  \underline{n}_i^j\nabla\ell_j(\underline{y}_i^j)^\intercal  = \underline{n}_i^j Q S\nabla'\gamma_j(T_j(\underline{y}_i^j))^\intercal   
  = 1, \quad \underline{y}_i^j\in Y_j^\rho,
  \label{eq:geomcond2}
\end{equation}  
where $\underline{n}_i^j$ is the initial normal at $\underline{y}_i^j$.

\begin{figure}[!htb]
  \centering
  \begin{tikzpicture}
    \node at (0,0) {\includegraphics[width=0.5\textwidth]{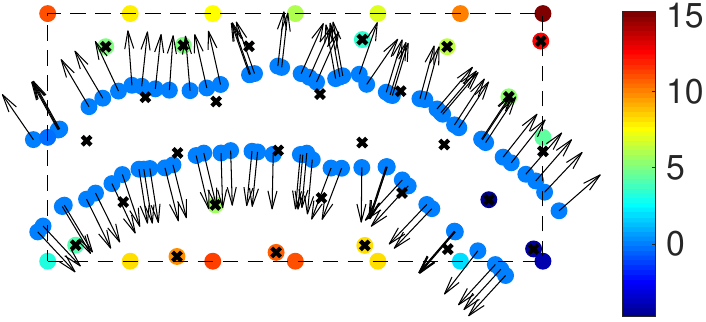}};
    \node at (-0.5,-1.3) {\large $2R_0$};
    \node at (-3.5,0.3) {\large $H_0$};
  \end{tikzpicture}  
  \caption{The local least-squares problem for one patch in the approximation of the 2D diaphragm data. The center points ($\times$), the data points (colored dot and arrow), and the additional points $X_{j,e}$ (colored dot) are shown. The color corresponds to the function value and the arrows are the normal directions used in conditions~\cref{eq:geomcond2} and~\cref{eq:geomcond3}.}
  \label{fig:extremal}
\end{figure}

The conditions~\eqref{eq:geomcond1} and~\eqref{eq:geomcond2} are all positioned at the object surface. This implies that in regions of the patch that are far from the object, the level set function may misbehave, and for example introduce spurious zero level set curves. To stabilize the approximation, we add a third set of conditions. We start from the set $X^\prime$ of RBF centers in the reference patch. We compute the average node distance $h^\prime=\frac{1}{n}\sum_{i=1}^n \min_{j\neq i}\|x_i^\prime-x_j^\prime\|$ (which we use as an approximation of the fill distance) and use this to construct a node set $X_b^\prime$ with approximate node distance $h'$ at the reference patch surface. We combine the boundary nodes with centers that are at least a distance $\alpha h^\prime$ away from the object to get the set 
$X_{j,e}^\prime=X_b^\prime\cup\left\{\underline{x}_k^\prime\in X^\prime\,|\,\min_i\|\underline{x}_k^\prime - T_j(\underline{y}_i^j)\|>\alpha h^\prime\right\}$. Then, for each $\underline{y}^{j,e}_i$ such that $T_j(\underline{y}^{j,e}_i)\in X^\prime_{j,e}$ we find the index of the closest data point on the object $\pi_{i,j} = \argmin_k\|\underline{y}^{j,e}_i-\underline{y}_k\|$, where $\underline{y}_k\in Y$. We use the approximate distance from the boundary data in the normal direction as a condition for the value of the level set function at that point
\begin{equation}
  \ell_j\left(\underline{y}^{j,e}_i\right) = \gamma_j(T_j(\underline{y}^{j,e}_i))= \left(\underline{y}^{j,e}_i-\underline{y}_{\pi_{i,j}}\right)\cdot \underline{n}_{\pi_{i,j}}, \quad T_j(\underline{y}^{j,e}_i)\in X^\prime_{j,e}.
  \label{eq:geomcond3}
\end{equation} 
The three types of conditions are illustrated for a two-dimensional patch in~\cref{fig:extremal}.
Equations~\eqref{eq:geomcond1}, \eqref{eq:geomcond2}, and~\eqref{eq:geomcond3} form the local least-squares problems. Solving these gives the vectors $\underline{u}_j$ of patch nodal values, which allows us to evaluate the local RBF approximation~\eqref{eq:RBF}. Using~\eqref{eq:uofv} and~\eqref{eq:rbfpum} then provides the global level set function at all points within the union of patches. 


The algorithm that we have derived has a number of parameters that we discuss qualitatively here, and more quantitatively in~\cref{sec:exp}. The number of patches $P$ together with the relative overlap $\delta_0$ determine the number of data points $m_j$ within a patch. Furthermore, $\delta_0$ governs the transition between the solutions in adjacent patches. A small $\delta_0$ can lead to large gradients in the transition regions, while a large $\delta_0$ increases the computational cost. The number of center points $n$ in the reference patch and the overlap $\delta_0$, together with the number of patches $P$, determine the oversampling, and accordingly how closely the local data is fitted.    

The reference patch dimensions $R_0$ and $H_0$ are coupled parameters. Therefore, we use a scale parameter $R_s$, such that $R_0=R_s\bar{R}$ and an anisotropy parameter $H_s$, such that $H_0=H_sR_s\bar{H}$. When performing numerical experiments, we fix the scale parameter, and only vary the anisotropy parameter.


The multiquadric RBF $\phi(r)=\sqrt{1+r^2}$ can be equipped with a shape parameter $\varepsilon$ by using $\phi(\varepsilon r)$ as the basis functions. A small shape parameter is appropriate when approximating a function that is smooth in the sense that it mainly/only contains low frequency Fourier modes. At the same time, the conditioning of the local problem grows with an inverse power of $\varepsilon$~\cite{Schaback05,LaFo05}. The ill-conditioning issues can be removed by using a stable evaluation method~\cite{FoWri04,FoPi07,FoLaFly11,FassMc12,LLHF13,FoWri17,KoLaYu19}. However, here we choose to use the multiqadric RBF with an intermediate shape parameter value that gives a reasonable condition number. This makes the local approximations robust and well-behaved also away from the object surface. 

The parameters $\rho$ and $\alpha$ that determine which extra conditions to use outside and inside the domain are less crucial for the results, but the existence of extra conditions is important for the robustness of the local problem.

\begin{figure}[!htb]
  \centering
  \begin{tikzpicture}
  \node at (3,0)  {\includegraphics[width=0.7\textwidth]{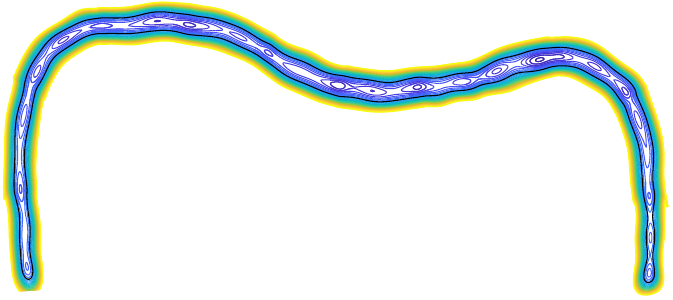}};
  \node at (2,-1) {\includegraphics[width=0.7\textwidth]{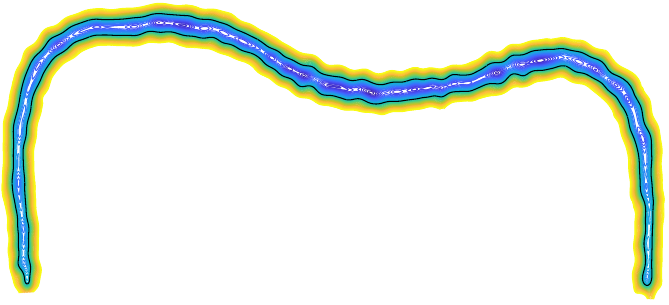}};
  \node at (1,-2) {\includegraphics[width=0.7\textwidth]{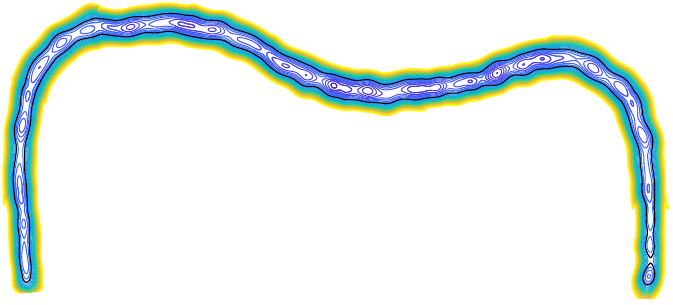}};
  \node at (0,-3) {\includegraphics[width=0.7\textwidth]{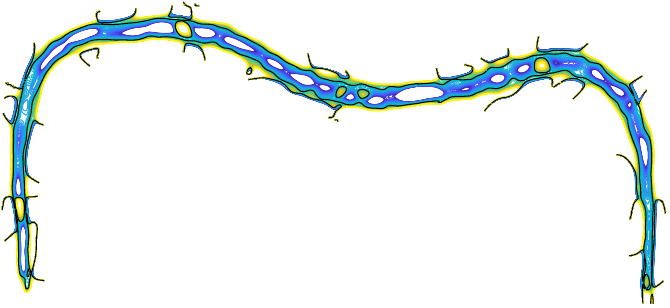}};
  \end{tikzpicture}
\caption{Level set functions for the 2D diaphragm for different parameter choices. The zero level set, i.e., the surface of the reconstructed geometry is the black curve. The result for $H_0=1$, $R_0=1$, $n=21$, $P=12$, $\delta_0=0.25$, $\varepsilon=1$, $\rho=0.25$, $\alpha=0.5$ is shown at the back. From back to front, for the other cases, the differences are for the first case $n=91$, for the second case $H_0=0.25$, and for the third case $n=91$ and $X_{j,e}=\emptyset$.}
\label{fig:geomparam}
\end{figure}

The effect of the most important parameters are shown in~\cref{fig:geomparam}. We see that increasing $n$ leads to a closer fit of the data, which in this case means fitting the noise. When $H_s<1$, it becomes hard to differentiate between the two surfaces, and the algorithm makes a mistake at the lower right part and creates a disconnected zero level set. When we remove the extra conditions inside the patch and use a large $n$ the resulting level set function has spurious zero level sets both inside and outside the actual geometry. This effect is almost invisible for $n=21$, but worsens with increasing $n$. In these experiments, the local node set $X^\prime$ consisted of Halton nodes~\cite{Halton60}.

\subsection{Geometry quality measures}\label{sec:qualmeas}
To evaluate the reconstructed geometries numerically, we introduce two different quality measures. As a measure of smoothness, we use the mean Frobenius norm of the surface Hessian of the geometry, which is computed using the following steps:
\begin{enumerate}
\item The data points $\{\underline{y}_k\}_{k=1}^M$ are iteratively moved onto the geometry surface using the values and gradients of the levelset function $\ell(\underline{y})$.
\item The Hessian $\mathcal{H}_\ell(\underline{y})$ is computed at each surface point $\{\underline{y}^S_k\}_{k=1}^M$ using~\eqref{eq:transder}.
\item For each surface point, we find a rotation $Q_y$, such that the normal at $y$ becomes the $d$th unit vector. The surface Hessian is then given by
  \begin{equation*}
    \mathcal{H}^S_\ell(\underline{y})=\left(Q_y\mathcal{H}_\ell(\underline{y}) Q_y^\intercal\right)_{1:d-1,1:d-1}.
  \end{equation*}
\item Finally, the mean Frobenius norm, $F_H$ is computed as
  \begin{equation*}
    F_H = \frac{1}{M}\sum_{k=1}^M\|\mathcal{H}^S_\ell(\underline{y}_k^S)\|_F.
    \label{eq:frobNorm}
  \end{equation*}  
\end{enumerate}
When least-squares approximation is employed for smoothing, the resulting smoothness is coupled to the amount of oversampling. That is, we expect $F_H$ to decrease if we decrease the number of degrees of freedom.  
However, there is also a problem dependent lower bound for how many degrees of freedom we need, to capture the true structure behind the noisy data. The second quality measure is an indicator specific to thin structures, which we have found to be effective in detecting when issues related to lack of resolution appear. It requires the nodes to be annotated with inner, outer and edge, which can be done automatically. We let $K_\mathrm{in}$ be the index set corresponding to data points on the inner surface, and we let $\tau(\underline{y}_k)$ be the distance to the closest data point on the outer surface. We measure the relative loss of thickness $T_L$ as  
\begin{equation}
  T_L = 1-\min_{k\in K_\mathrm{in}}\frac{\tau(\underline{y}^S_k)}{\tau(\underline{y}_k)}.
  \label{eq:thickMeasure}
\end{equation}
The loss is not expected to become zero, but in a well-resolved situation it is on the order of the noise size.

\section{Numerical results for the geometry reconstruction}\label{sec:exp}
The implementation is written in MATLAB and all experiments have been performed on standard laptop computers.

In the experiments we have chosen to use Halton nodes~\cite{Halton60} for the RBF center points $X'$ in the reference patch. The main reason is that low discrepancy points are insensitive to scaling and are therefore well-distributed both in the reference patch and in the physical domain. In the experiments, we have fixed $R_s$ such that $R_0=1$ for $P=10$ in 2D, and for $P=20$ in 3D.

The set of parameters we can use to optimize the quality of the reconstructed geometry are listed and described in \cref{tab:1.5}

\begin{table}[!htb]
  \centering
  \footnotesize
  \caption{The tunable parameters and their descriptions}
  \label{tab:1.5}
  \begin{tabular}{|c|p{0.8\textwidth}|}\hline
    $\alpha$ & Determines the mimimum distance from the surface to the extra points in~\eqref{eq:geomcond3}. Increasing $\alpha$ will slightly decrease the oversampling factor. \\\hline
    $\rho$ & decides how far outside the patch boundaries data is collected for the local problems to form the data conditions~\eqref{eq:geomcond1} and~\eqref{eq:geomcond2}. Increasing $\rho$ increases the oversampling slightly, but the main and desired effect is to stabilize the behavior of the level set function where the patch intersects the geometry.\\\hline
    $H_s$ & Controls the aspect ratio of the reference patch. By increasing $H_s$ we make the contribution to the residual from the distance normal to the surface larger than the contribution from the distance tangential to the surface. That is, we increase the amount of surface smoothing. Also, the distance between the two surfaces as well as the end curvature is decreased in the computational domain.\\\hline
    $n$ & This is the number of RBF center points in each patch, which has the strongest influence on the total number of degrees of freedom.  \\\hline
    $P$ & This is the number of patches, which is the second most important parameter in terms of degrees of freedom.\\\hline
    $\delta_0$ & is the relative overlap between patches. When $\delta_0=1$ patch neighbors on opposite sides of one patch meet in the middle. Increasing $\delta_0$ increases the oversampling, since more data points are included in the local subproblems, while the number of RBF center points in the patch is the same. The resolution of the geometry is however decreased, since the same number of basis functions cover a larger area in the physical domain.\\\hline
    $\varepsilon$ & The RBF shape parameter affects the surface Frobenius norm. With a larger shape parameter more variations in the data, including the noise, are picked up. \\\hline
  \end{tabular}  
\end{table}

\subsection{Results for 2D cross sections}
For the 2D cross section data shown in~\cref{fig:data} one parameter at a time was varied to study its behavior for all combinations of $n=21$, $28$, and $36$ and $P=10$, $20$, $30$, and $40$. In this phase, we mainly focus on the smoothness measure $F_H$.

The smoothness of the reconstruction is relatively insensitive to the distance $\rho h$ at which we pick up data outside the patch. The choice $\rho=0.5$ performs well in all cases, and is significantly better than $\rho<0.5$, especially for larger numbers of patches. We keep this parameter fixed in all subsequent experiments.
The reconstruction is even less sensitive to the distance $\alpha h$ between the object and the extra conditions inside the patch. We fix this parameter to $\alpha=1.5$.
For the overlap, the smoothness increases significantly for the range from $\delta_0=0.1$ to $0.9$. However it is worth noticing that the loss in thickness increases with $\delta_0$. Depending on the data, different overlaps may be optimal. In the later experiments, we use both $\delta_0=0.5$ and $\delta_0=0.9$.
The curvature measure $F_H$ decreases with the shape parameter down to $\varepsilon\approx 0.2$ for the tested cases. For small values the thickness loss increases and the computations gradually become ill-conditioned, since we are solving the local systems directly. To stay on the robust side, we use $\varepsilon=0.5$ in the subsequent experiments.
The patch scaling is fixed to $R_s\approx 0.0215$, which makes $R_0=1$ for $P=10$ patches. An increase of the anisotropy factor $H_s$ leads to a smoother result, especially for larger $P$, whereas loss of thickness occurs to a larger degree for small $H_s$. In the later experiments, we use both $H_s=1$ and $H_s=2$.






To find reconstruction parameters that minimize the quality measures we optimized $F_H+T_L$ over $P_0$ for different $n$ and the four combinations of overlap and aspect ratios. Since $T_L$ is an order of magnitude larger than $F_H$, the priority becomes to find a reconstruction that keeps an even width throughout. The resulting parameters are reported in~\cref{tab:1}. We note that all choices of $n$ gave good reconstruction results. If the even smaller $n=6$ is used, we also get a good result in terms of quality, but points move more with respect to the given data. 
\begin{table}[!htb]
  \caption{The best results for different $n$ if $F_H+T_L$ is minimized over $P_0$ for the four combinations of $\delta_0$ and $H_s$.}\label{tab:1}
  \footnotesize
  \centering
\begin{tabular}{r|rrrrrr}\hline
  $n$ & $P_0$ & $P$ & $\delta_0$ & $H_s$ & $F_H$ & $T_L$\\\hline
  10 & 30& 34 & 0.5 & 1 & 0.024 & 0.20\\
  15 & 33& 37 & 0.5 & 1 & 0.022 & 0.14\\
  21 & 22& 24 & 0.5 & 1 & 0.016 & 0.11\\ 
  28 & 12& 12 & 0.5 & 1 & 0.013 & 0.13\\
  36 & 11& 11 & 0.5 & 1 & 0.014 & 0.12\\
  45 & 11& 11 & 0.5 & 1 & 0.014 & 0.12\\ \hline
\end{tabular}
\end{table}

\begin{figure}
  \begin{tabular}{c@{}}
  \includegraphics[width=0.44\textwidth]{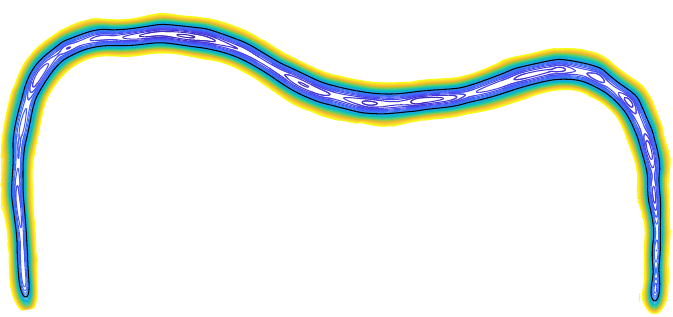}\hspace*{1cm}\\
  \includegraphics[width=0.49\textwidth]{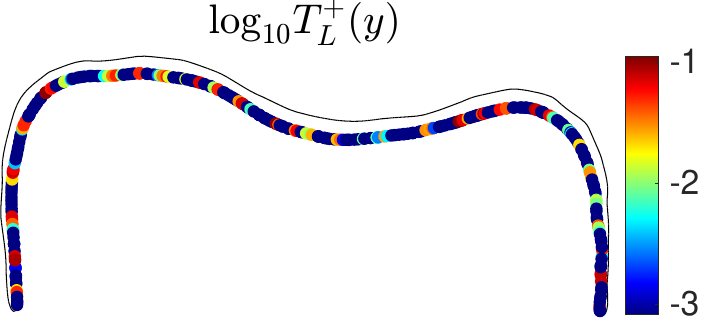}\\
  \includegraphics[width=0.49\textwidth]{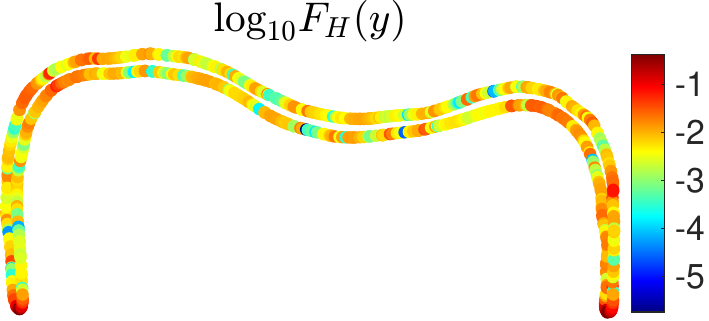}
  \end{tabular}%
  \begin{tabular}{@{}c}
  \hspace*{4mm}\includegraphics[width=0.42\textwidth]{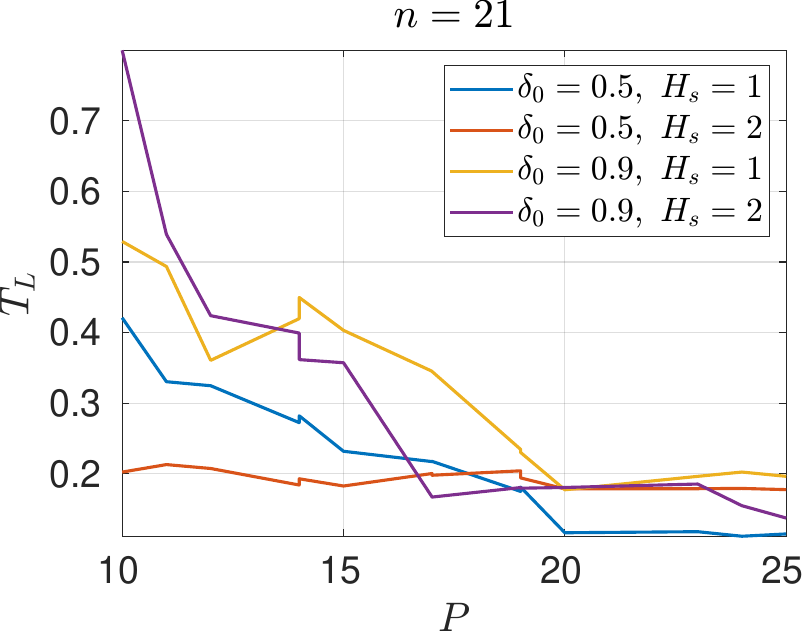}\\
  \includegraphics[width=0.45\textwidth]{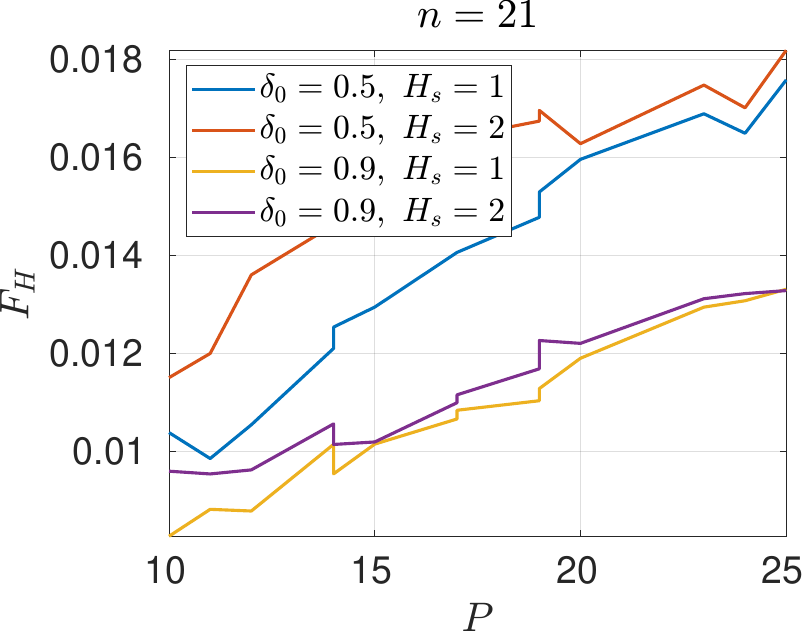}
  \end{tabular}
  \caption{The best reconstruction for $n=21$ and $P=24$ with $F_H=0.016$ and $L_T=0.11$ is shown together with the spatial distribution of the quality measures to the left. The black curve is the zero level set. Negative losses are set to zero and correspond to dark blue. The quality measures are shown as a function of the number of patches $P$ to the right.}\label{fig:2Dqual}
\end{figure}

The best result according to the objective function is for $n=21$. The geometry and the quality measures are illustrated in \cref{fig:2Dqual}. In the left side of the domain, the loss of thickness seems to be correlated with the patch overlap regions, while this is less clear in the rest of the domain. In the right part of the figure, we can observe that if we allowed for a higher thickness loss, we could reach an overall smoother result.

\begin{figure}
  \centering
  \includegraphics[width=0.40\textwidth]{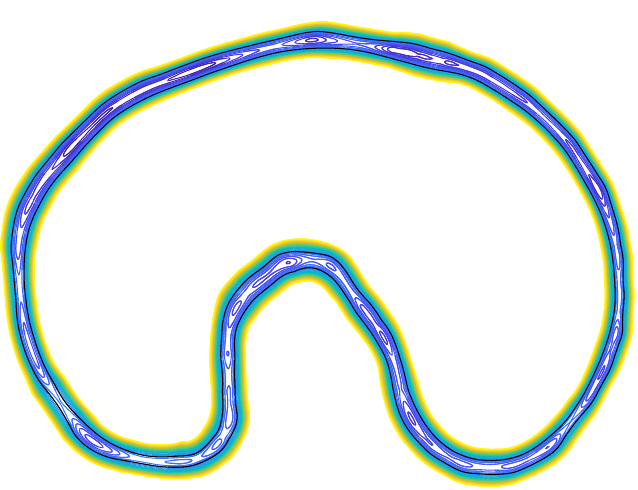}\hspace*{0.07\textwidth}
  \includegraphics[width=0.40\textwidth]{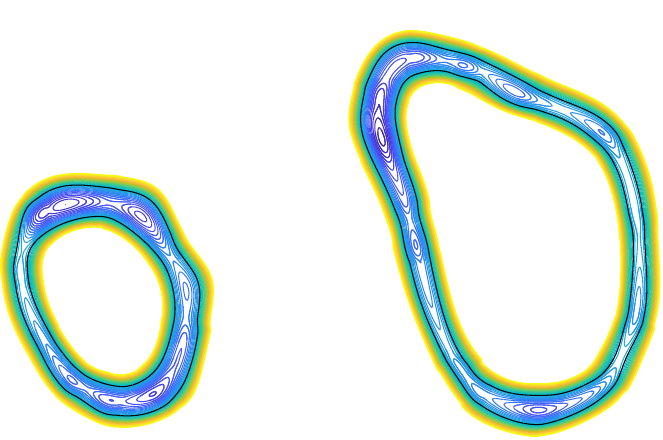}
  \caption{Reconstructed diaphragm cross sections for $n=21$ in two other planes. The best parameters for the left case are $P_0=40$, $P=46$, $\delta_0=0.9$, and $H_s=1$, leading to $F_H=0.008$ and $T_L=0.22$. The best parameters for the right case with two loops are $P_0=21$, $P=23$, $\delta_0=0.5$, and $H_s=1$, leading to $F_H=0.013$ and $T_L=0.14$.}\label{fig:2Dother}
\end{figure}

To show that the reconstruction algorithm can handle different geometries, we run the same code with the same objective function for two other cross sections with $n=21$. The results including best parameters and quality measures are shown in \cref{fig:2Dother}.


\subsection{Results for 3D reconstruction}
\begin{figure}[!htb]
  \centering
  \includegraphics[width=0.40\textwidth]{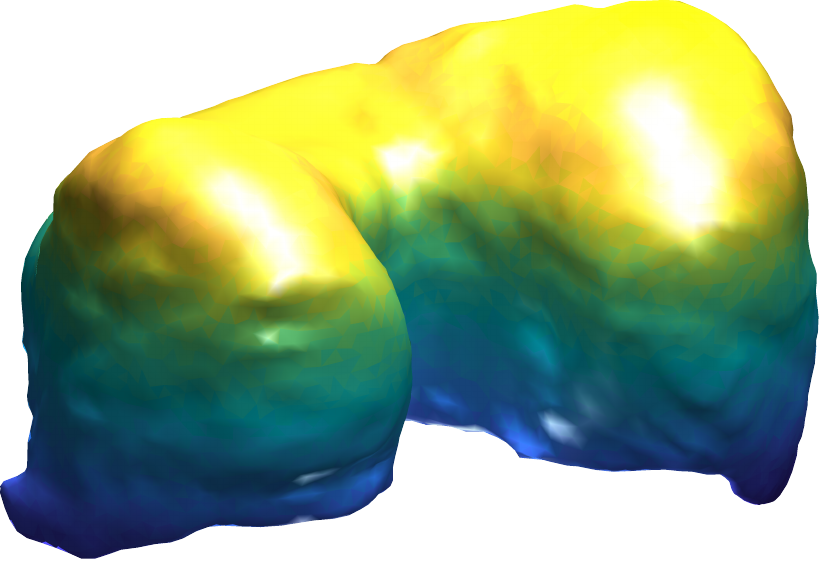}\hspace*{0.07\textwidth}
  \includegraphics[width=0.40\textwidth]{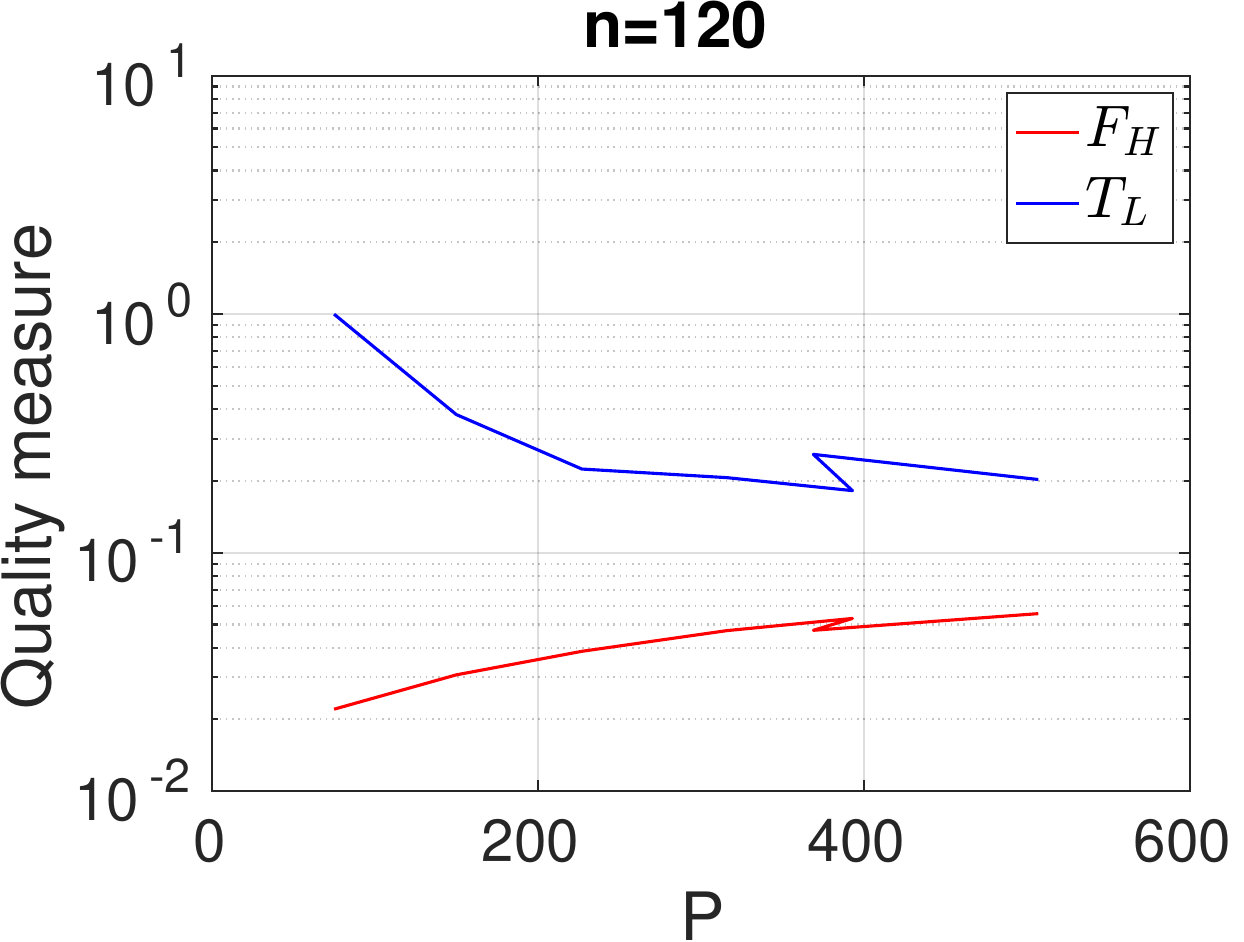}
  \caption{The best 3D diaphragm reconstruction, using $n=120$ and $P = 393$ leading to $F_H = 0.053$ and $T_L = 0.18$ (left). The quality measures are shown as a function of the number of patches, with $n = 120$ (right).}\label{fig:3DProb3}
\end{figure}

\begin{figure}[!htb]
  \centering
  \includegraphics[width=0.40\textwidth]{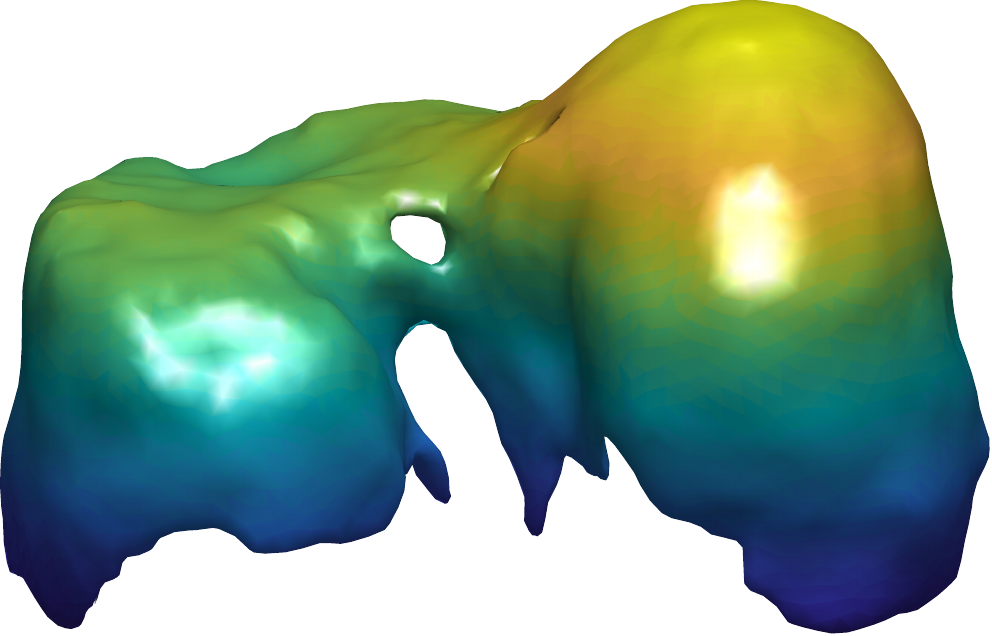}\hspace*{0.07\textwidth}
  \includegraphics[width=0.45\textwidth]{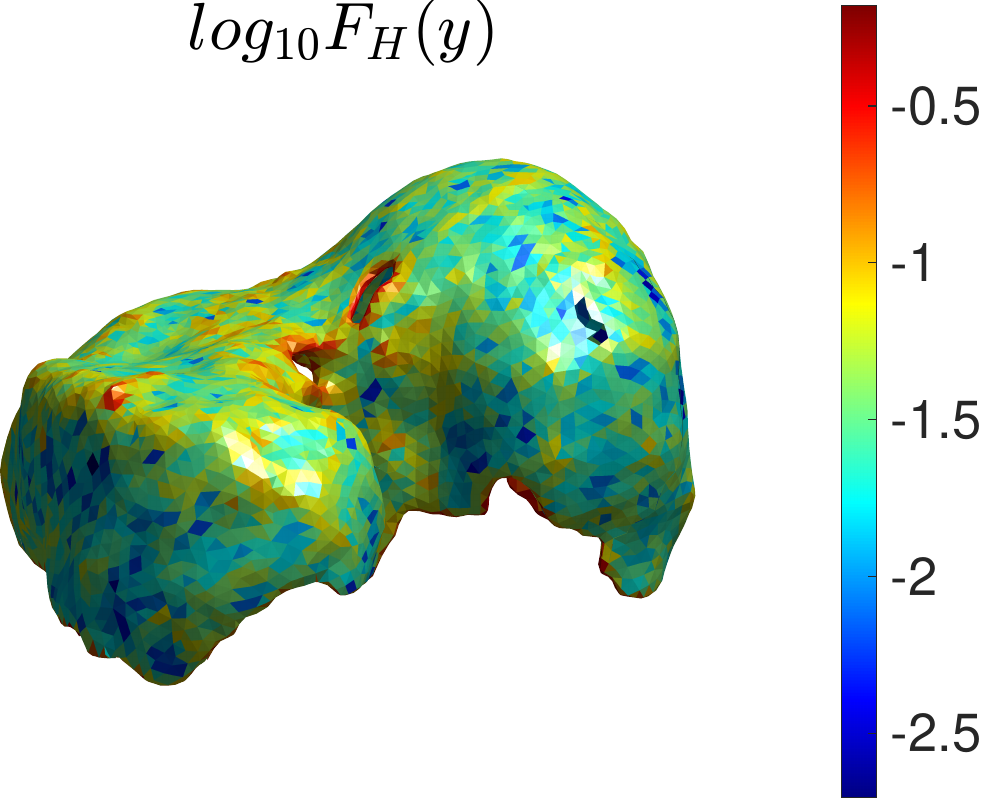}
  \caption{The second 3D diaphragm geometry reconstruction for $n=120$ and $P = 393$ with $F_H = 0.067$ and $T_L = 0.72$ (left). The Frobenius norm plotted on the reconstructed surface of the diaphragm (left).}\label{fig:3DProb5}
\end{figure}

We reconstruct 3D diaphragm geometries from two different patients. The reconstructions are shown in \cref{fig:3DProb3} and \cref{fig:3DProb5}. Similarly to the reconstruction of the 2D cross sections we optimize for each parameter separately while fixing all others. We choose parameters that minimize the quality measure $F_H + T_L$ and all parameter choices are made using the geometry shown in Figure \cref{fig:3DProb3}. The same parameters are used to reconstruct the more complex geometry.

The reconstruction is quite robust in terms of the distance for which we pick points outside the patch, $\rho h$ and the distance $\alpha h$ between the object and the extra conditions inside the patch. Here, the quality measure is minimized when $\alpha = 0$ and $\rho = 0$. We additionally find that the quality measure is fairly constant with respect to the chosen anisotropy ratio $H_s$ of the reference patch, with slightly better results provided when $H_s = 1$. Moreover, similar to the 2D reconstruction, an increase in relative overlap $\delta_0$ results in a significant increase in smoothness. However, it also results in a loss of thickness and the optimal value based on our metric is $\delta_0 = 0.4$. The optimal shape parameter $\varepsilon$ with respect to the resulting smoothness is relatively small ($\varepsilon = 0.25$), however the loss of thickness is at a maximum ($T_L = 1$). This means that some data points originally annotated with being on the inner surface have moved to the outer surface or vise versa. The optimal quality measure is obtained when $\varepsilon = 0.9$.

We note that point flipping between inner and outer surfaces is unique to the 3D reconstruction problems, were the distance between the inner and outer surface data points, as well as the inner and outer surface reconstruction, can be very small in specific areas. This can indicate different issues, either relating to the annotation of points, the data quality, or the reconstruction quality. Hence, in cases when this happens at isolated points, which indicates bad data quality or annotation, we do not include these points in the thickness loss measure $T_L$.
However, in cases where multiple points in the same neighborhood exhibit the same issue, they cannot be ignored.

With the aforementioned parameters fixed, we vary the degrees of freedom to see how the quality measures are affected. We find the optimal number of patches $P_0$ for each choice of RBF center points $n$. Results are provided in~\cref{tab:2}. 
\begin{table}[!htb]
  \caption{The best results for different $n$ if $F_H+T_L$ is minimized over $P_0$ for fixed values of $\alpha=0$, $\rho=0$, $\delta_0=0.4$, $H_s=1$ and $\varepsilon=0.9$.}\label{tab:2}
  \footnotesize
  \centering
\begin{tabular}{r|rrrr}\hline
  $n$ & $P_0$ & $P$ & $F_H$ & $T_L$\\\hline
  35 & 300 & 393 & 0.031 & 0.32\\
  56 & 420 & 507 & 0.039 & 0.28\\
  84 & 300 & 393 & 0.044 & 0.21\\ 
  120 & 300 & 393 & 0.053 & 0.18\\
  165 & 180 & 227 & 0.046 & 0.21\\
  220 & 180 & 227 & 0.050 & 0.19\\ 
  \hline
\end{tabular}
\end{table}

We obtain the optimal result according to the objective function when $n=120$. The two reconstructed diaphragm geometries are shown in \cref{fig:3DProb3} and \cref{fig:3DProb5}. In \cref{fig:3DProb3} it is also shown how the two quality measures change when the patch number is increased. As is expected from the 2D examples, an increase in the number of patches $P$ results in an increase in the Frobenius norm $F_H$, but a decrease in thickness loss $T_L$. For the more complex diaphragm geometry we also provide the smoothness measure over the surface. Note the high $F_H$ values around the holes, where the vena cava and esophagus pass through the diaphragm, which make this a more challenging geometry to reconstruct.

\section{The PDE solution}\label{sec:PDE}
Given the smooth geometry reconstruction, we want to use it for solving a PDE problem. In particular, we eventually want to model the biomechanical action of the diaphragm during respiration, see also~\cite{invive,TVLBC21}. The base method, RBF-PUM, is the same as for the reconstruction, but there are some differences which we describe in detail in this section. In this paper we solve a Poisson problem with Dirichlet boundary conditions as a proof of concept. We manufacture the right hand side functions based on an infinitely smooth solution function $u$ and write the problem as
\begin{eqnarray}
\Delta u(\underline{x}) &=& f_i(\underline{x}),\quad \underline{x}\in\Omega,\label{eq:pde} \\
u(\underline{x}) &=& f_b(\underline{x}),\quad  \underline{x}\in\partial\Omega.\label{eq:bc}
\end{eqnarray}
%
We sample the PDE and boundary condition at the point sets $Y_i$ and $Y_b$, respectively. To form the global linear system, we use~\cref{eq:PDERBFPUM} with $\mathcal{L}_k=\Delta$ and $f=f_i$ for $\underline{y}_k\in Y_i$ and $\mathcal{L}_k=I$ and $f=f_b$ for $\underline{y}_k\in Y_b$. For the derivatives of the local approximations we use \eqref{eq:transder} and for the derivatives of the weight function we also use~\cref{eq:weightder}. The scaling of the equations is implemented as described in~\cref{sec:rbfpum}, using the power $h^{-1/2}$ for the Dirichlet conditions.

%
%
To generate a cover for the PDE computation, we need a point set on the computed object surface. We use the initial data point cloud from the reconstruction problem, and move these points iteratively onto the surface using the gradient of the level set function. This works well when the level set function is free of spurious zero level set curves. Then we run the cover algorithm derived in~\cref{sec:cover}. 

To form the interior node set $Y_i$, we first generate Halton nodes within each patch. Then we use the level set function to determine which points are inside and which are outside, and finally we remove points that are close to each other. Least-squares RBF methods are not sensitive to the exact distribution of the sampling points as long as these are quasi uniform~\cite{LaShchHer17,ToLaHe21,LaMaMiPoo22}.
For the surface point set $Y_b$, we use the mesh representation of the initial point cloud. We refine or decimate the mesh and move all new points onto the surface. It would be possible to construct a surface node generation algorithm that uses the level set function directly, but we have not investigated this aspect here.

Using RBF-PUM, there are two standard modes of refinement. Patch refinement with a fixed number of nodes per patch leads to algebraic convergence and node refinement with fixed patches leads to spectral convergence for infinitely smooth RBFs~\cite{LaShchHer17}. For both refinement modes, the shape parameter should be kept fixed. If a stable evaluation method is not used, the ill-conditioning of the local subproblems increases with the refinement level~\cite{LaShchHer17}. Therefore, we use Gaussian RBFs and the RBF-QR method~\cite{FoLaFly11,LLHF13} for stable evaluation when solving PDEs. We note that to derive the transformed Laplacian~\cref{eq:transder}, we need the second derivatives with respect to each coordinate direction. These derivatives were not present in the 3D RBF-QR implementation derived in~\cite{LLHF13}, but we have now upgraded the implementation with all second derivatives~\cite{RBFQR3Dmixed}. Furthermore, we replace the infinitely smooth bump function with the $C^2$ Wendland function~\cref{eq:wend} as generating function for the weights. This implies that we only require $C^2$ continuity of the global solution, while retaining the global high-order convergence arising from the high-order local approximations~\cite{Wend02,Fass07,LaShchHer17}. 

In the convergence tests performed here we use the patch refinement approach which is the most appropriate for large scale problems. It was shown numerically in~\cite{LaShchHer17} that least-squares RBF-PUM is stable under patch refinement.
For our particular setting, we need to keep in mind that the computational domain is (several instances of) the reference patch. By keeping the scale factors $R_s$ and $H_s$ constant under refinement, we ensure that refining the physical domain corresponds to refining the computational domain.
We choose to do the refinment only in the radial direction. That is, we keep the structure with one layer of patches.
The node layout in the reference patch as well as the physical domain is quasi uniform, and as a practical measure of the worst case fill distance in the physical domain, we use
\begin{equation}
  h = \max_{1\leq j\leq P}\left(\frac{|\Omega_j|}{n}\right)^\frac{1}{d}=C_d\max_{1\leq j\leq P}\left(\frac{H_jR_j^{d-1}}{n}\right)^\frac{1}{d}.
  \label{eq:hell}
\end{equation}
Since $H_j$ should be proportional to the width of the geometry and almost constant, the fill distance is approximately proportional to $R_j^{(d-1)/d}$. That is, by dividing a patch in half in 2D, we only reduce the fill distance by $1/\sqrt{2}$.
%

Let $\mathcal{I}_h(u)=\sum_{j=1}^Pw_j\mathcal{I}_h(u_j)$, where $u_j=u|_{\Omega_j}$, be the RBF-PUM interpolant of the function $u$, and note that $u(\underline{x})=\sum_{j=1}^Pw_ju_j$, $\underline{x}\in\Omega$ due to the partition of unity property. Then we can express the interpolation error as 
\begin{equation}
  \mathcal{E}_I=\mathcal{I}_h(u)-u = \sum_{j=1}^Pw_j(I_h(u_j)-u_j)\equiv \sum_{j=1}^Pw_je_j,\quad \underline{x}\in\Omega.
  \label{eq:ei}
\end{equation}
%
%
Similarly, the consistency error for the Laplacian operator becomes~\cite{Wend02,LaShchHer17}
\begin{equation}
  \mathcal{E}_\Delta=\Delta\sum_{j=1}^Pw_je_j =\sum_{j=1}^P(\Delta w_je_j+2\nabla w_j\cdot\nabla e_j+w_j\Delta e_j) ,\quad \underline{x}\in\Omega.
  \label{eq:deltae}
\end{equation}
The derivatives of the $C^k$ partition of unity weight functions are typically expressed in terms of the radius of the patch, especially in a setting with spherical patches~\cite{Wend02,LaShchHer17}. In our context with a non-uniform distribution of cylindrical patches, we need to be more careful. The largest size of the derivatives depends on the smallest absolute overlap $\delta$ between patches. We use the following estimate
\begin{equation}
  \|D^\beta w\|_\infty =G_{|\beta|}\delta^{-\beta},\quad |\beta|\leq k.
  \label{eq:west}
\end{equation}
%
%
For the local interpolation error, we use results derived in~\cite{RieZwi10}.
For our specific setting, we first make the estimate in the computational domain, and in a second step include the transformation. We let $v_j(x^\prime)=u_j(x)$ and define $n_p$ as the dimension of a polynomial space of degree $p$ in $d$ dimensions. Then we use the following error estimate;
\begin{equation}
  \|D^\beta (I_h(v_j)-v_j)\|_\infty \leq C_{|\beta|} h^{p+1-|\beta|-\frac{d}{2}}\|v_j\|_{\mathcal{N}(\Omega_0)},
    \label{eq:eest}
\end{equation}
where $n_p\leq n <n_{p+1}$, the constant $C_{|\beta|}$ depends on the shape of the reference patch and on $d$, and $\|\cdot\|_{\mathcal{N}(\cdot)}$ denotes the native space of the infinitely smooth radial basis function that is used.
\begin{remark}From the theoretical perspective, the results in~\cite{RieZwi10} require $h$ to be small enough for the estimates to be valid. In our setting $h$ is related to $n$ and cannot be decreased independently. The theoretical upper bound on $h$ is restrictive, while in numerical experiments the optimal rate is achieved also for larger $h$~\cite{LaShchHer17}.
\end{remark}
We finally use \eqref{eq:transder} to provide estimates for $e_j$ in the physical domain
\begin{equation}
  \|\nabla e_j\|_\infty\leq \|\nabla^\prime(I_h(v_j)-v_j)S^\intercal Q^\intercal\|_\infty\leq d\max_{1\leq i\leq d}s_{ii}C_1h^{p-\frac{d}{2}}\|v_j\|_{\mathcal{N}(\Omega_0)},
  \label{eq:egrad}
\end{equation}
\begin{equation}
  \|\Delta e_j\|_\infty \leq \sum_{i=1}^ds_{ii}^2 \|\partial^2/\partial x_i^{\prime 2}(I_h(v_j)-v_j)\|_\infty\leq d\max_{1\leq i\leq d}s_{ii}C_2h^{p-1-\frac{d}{2}}\|v_j\|_{\mathcal{N}(\Omega_0)}.
  \label{eq:edelta}
\end{equation}
%
%
Combining \cref{eq:ei}--\cref{eq:edelta} we get the consistency estimates
%
%
\begin{equation}
  \|\mathcal{E}_I\|_\infty\leq KG_0C_0h^{p+1-\frac{d}{2}}\max_j\|v_j\|_{\mathcal{N}(\Omega_0)},
\end{equation}  
\begin{eqnarray}
  \|\mathcal{E}_\Delta\|_\infty\leq K\hspace{-1.5em}&&\left(dG_2\delta^{-2}C_0h^{p+1-\frac{d}{2}} +
  2dG_1\delta^{-1}d\max_is_{ii}C_1h^{p-\frac{d}{2}}\right.\nonumber\\
  &&\quad+ \left.G_0d\max_is_{ii}C_2h^{p-1-\frac{d}{2}} \right)\max_j\|v_j\|_{\mathcal{N}(\Omega_0)},
  \label{eq:theorconv}
\end{eqnarray}
where $K$ is the maximum number of patches overlapping at any given point. The best convergence rate that can be achieved is $\mathcal{O}(h^{p-1-\frac{d}{2}})$. This requires the absolute overlap $\delta$ to go to zero with the same rate as $h$ or slower than $h$.
 We see from \cref{eq:hell} that $h\propto R_j^\frac{d-1}{d}$, which with $\delta=\delta_0R_j$ means that the relative overlap $\delta_0$ must increase as $1/\sqrt[d]{R_j}$ during refinement to achieve the optimal rate. This is an effect of the refinement mode we have chosen where patches are split only in $d-1$ directions. When patches are split uniformly, keeping $\delta_0$ constant is instead optimal. The scaling of the Dirichlet conditions with $h^{-1/2}$ is small enough to not affect the optimal rate.
%
%

 In the convergence experiments, 
 the relative maximum norm error is approximated as
\begin{equation}
  E = \frac{\|\tilde{u}-u\|_{\ell_\infty(Y)}}{\|u\|_{\ell_\infty(Y)}}.
\end{equation}  

\section{Numerical results for the PDE solution}\label{sec:pdeexp}

\begin{figure}[!htb]
  \centering
  \includegraphics[width=0.45\textwidth]{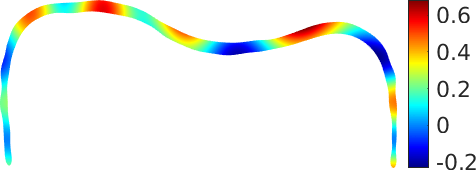}%
  \hspace*{1cm}\includegraphics[width=0.35\textwidth]{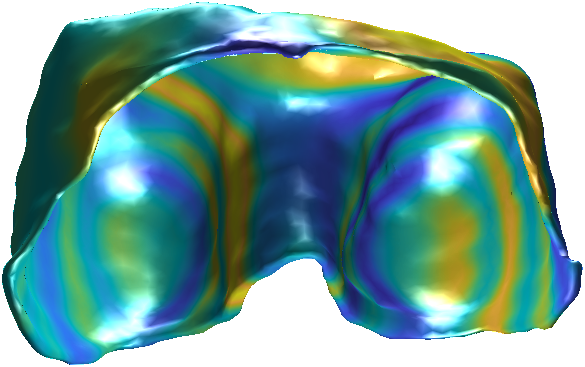}
  \caption{The manufactured solution functions $u_\mathrm{2D}$~\eqref{eq:u2D}
  and $u_\mathrm{3D}$~\eqref{eq:u3D}. For the latter, bright yellow corresponds to +1 and dark blue to -1.}
  \label{fig:sol2D}
\end{figure}

%
As the first test problem, we solve the Poisson problem~\eqref{eq:pde}--\eqref{eq:bc} in the 2D cross section of the diaphragm, with right hand side data generated from the solution function
\begin{equation}
  u_\mathrm{2D}=\frac{\sin(24(\tilde{x}_1-0.1)^2)\cos((\tilde{x}_2-0.3)^2) +\sin^2(6(\tilde{x}_2-0.5)^2)}{1+2\tilde{x}_1^2 + \tilde{x}_2^2},
  \label{eq:u2D}
\end{equation}
and in the 3D geometry with right hand side data generated from
\begin{equation}
  u_\mathrm{3D}=\frac{\sin(4\pi\tilde{x}_1\tilde{x}_2\tilde{x}_3) +\cos(4\pi(\tilde{x}_2-0.3)^2)}{1+2\tilde{x}_1^2 + \tilde{x}_2^2+0.5\tilde{x}_3^2}.
  \label{eq:u3D}
\end{equation}
These functions provide variation both along and across the diaphragm. The scaled variables are given by $\tilde{x}_i=x_i/\max_{x\in\Omega}|x_i|$. The solution functions are illustrated in \cref{fig:sol2D}. 

In the 2D experiments, which run up to a large number of patches, we have not increased the relative overlap as fast as suggested by the theory. The reason is that if the overlap becomes too large, patches that are not neighbours start to overlap, and this is an overlap that we cannot control, and which can produce large gradients in the weight functions. Instead of using $\delta_0\propto \bar{R}^{-1/d}$, we have used $\delta_0\propto \bar{R}^{-\frac{1}{2d}}$ for this case.

\begin{figure}[!htb]
  \centering
  \includegraphics[width=0.33\textwidth]{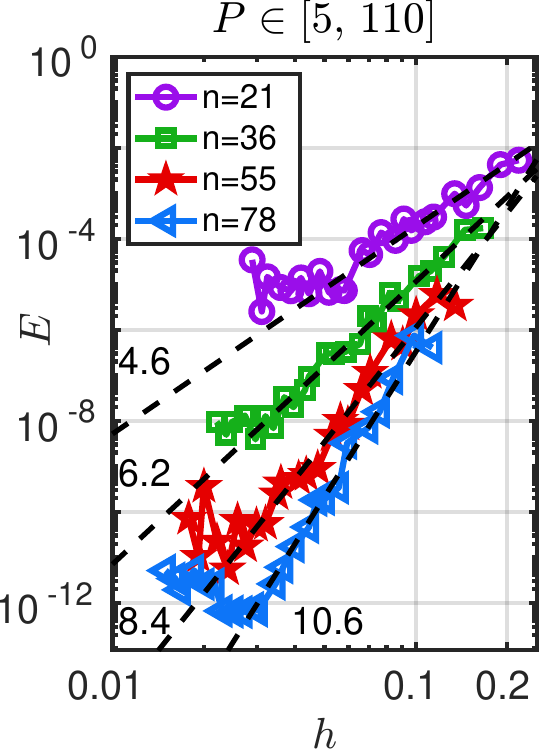}%
  \includegraphics[width=0.33\textwidth]{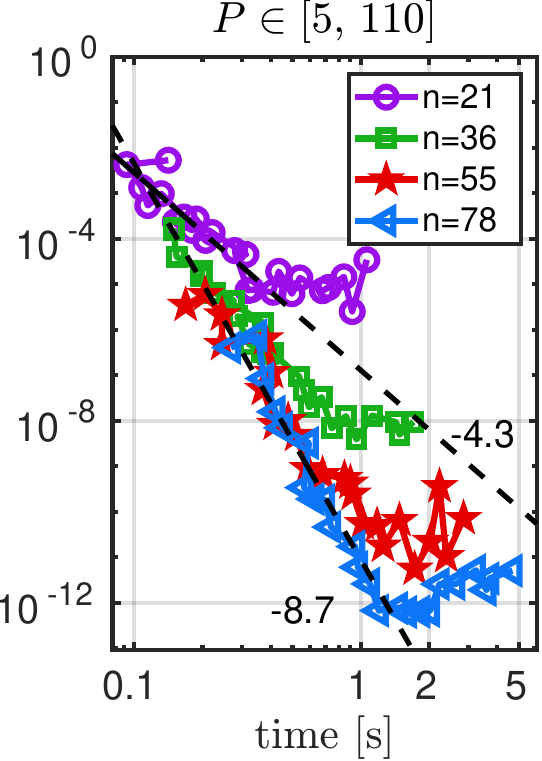}%
  \caption{The relative maximum error for the 2D problem is shown for different numbers of local points $n$ as a function of the fill distance $h$ (left subplot) and as a function of the run time in seconds (right subplot). The parameters not indicated in the subplots are chosen as the default values $q=3$, $\varepsilon=1$, $H_s=1$, and $\delta_0=0.1$. }
  \label{fig:Poisson2D}
\end{figure}

The purpose of the 2D experiment is to investigate how the result depends on the method parameters and how sensitive the parameter choices are. The left subplot of \cref{fig:Poisson2D} shows convergence as a function of the fill distance and the right one shows error as a function of time.
The best theoretical convergence rates from~\eqref{eq:theorconv} become 3, 5, 7, and 9 for these $n$. The numerical results seem to be around one order better.
We note that increasing the number of local points allows us to reach smaller errors for the same fill distance. To reach a smaller error in the most efficient way also requires gradually increasing $n$.

To investigate the parameter dependence, we have solved the PDE problem for all combinations of $n=21$, 28, 36, 45, 55, 78, and 91 and all 24 values of $P$ used for the plots for different parameter combinations. A measure of goodness was defined by looking for results that are more accurate and/or faster than the average performance when the other parameters are equal. Then for a selection of parameter values we compute in which percentage of cases each choice was the best and in what percentage of cases it failed, where fail means that the result is at least 5 orders of magnitude less accurate than for other parameter values.
\begin{table}[!htb]
  \centering
  \footnotesize
  \caption{The percentage of wins $B(\cdot)$, and fails $F(\cdot)$ for each tested parameter value. The shaded row indicates the default values used in the experiments.}
    \begin{tabular}{|cc|cc|ccc|ccc|}\hline
    $\varepsilon$ & $B(\varepsilon)$ & $H_s$ & $B(H_s)$ & $\delta_0$ & $B(\delta_0)$ & $F(\delta_0)$ & $q$ & $B(q)$ & $F(q)$\\\hline
    0.50  & 27 &      &    & & & &           2.0 & 51 & 5\\
    0.75 & 39 & 0.75 & 37 & 0.05 & 23 & 1 & 2.5 & 31 & 0\\
    \rowcolor{lightblue}
    1.00   &  34 &  1   & 51 & 0.1 & 61 & 0 &  3.0 & 18 & 0\\
        &    & 1.25  & 12 & 0.2 & 16 & 0 &      &    & \\\hline
  \end{tabular}  
\end{table}
Note that a smaller shape parameter value than 1 performs a little bit better, but the trend of the convergence curve was a bit more erratic, and we therefore picked $\varepsilon=1$ for improved stability. A smaller oversampling factor $q$ improves the efficiency, but the method can fail, so also here, we picked $q=3$ for stability. Also for the overlap, a small overlap is efficient, but when too small, the method can fail. It should also be noted that for all stable choices of parameters, we got a similar convergence and performance. That is, the main sensitivities can be understood and for the rest, the penalty for making a suboptimal choice is not too large.

\begin{figure}[!htb]
  \centering
  \includegraphics[width=0.49\textwidth]{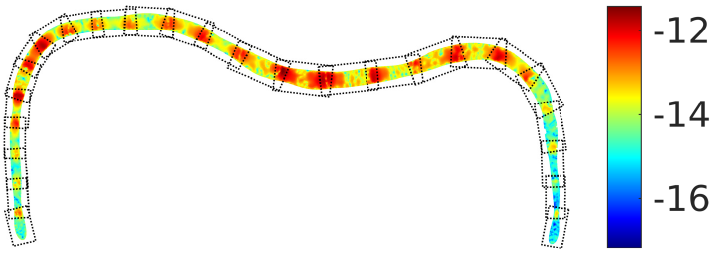}
  \includegraphics[width=0.49\textwidth]{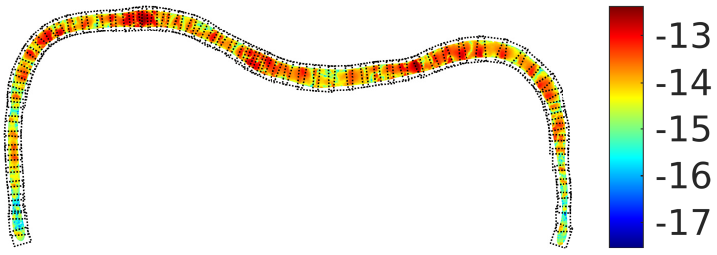}
  \caption{The logarithm in base 10 of the error in the RBF solutions for $n=78$ and $P=25$ (left) and $P=69$ (right).}
  \label{fig:Poisson2D_err}
\end{figure}

The error for two different numbers of patches is illustrated in \cref{fig:Poisson2D_err}. Especially in the left subfigure, it is clear that the error is a bit larger in the overlap regions, where the weight functions have large gradients.



\begin{figure}[!htb]
  \centering
  \includegraphics[height=0.38\textwidth]{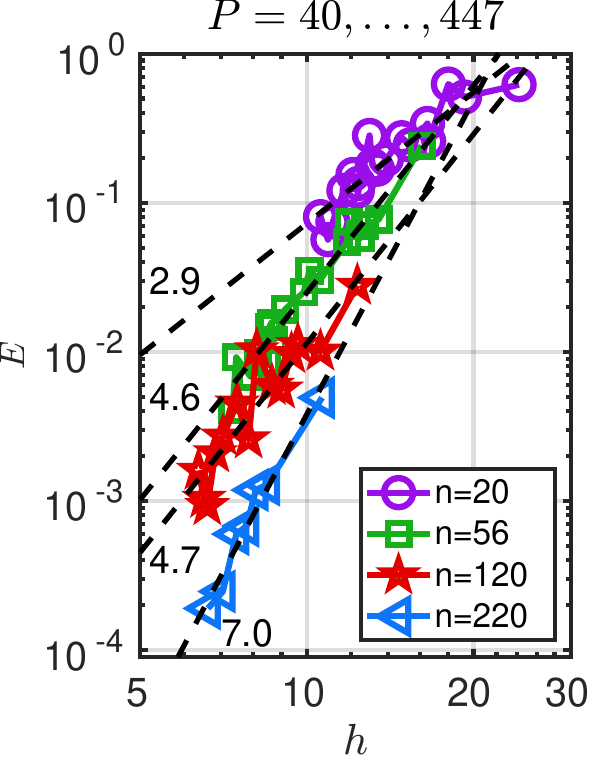}\hspace{1cm}
  \includegraphics[height=0.38\textwidth]{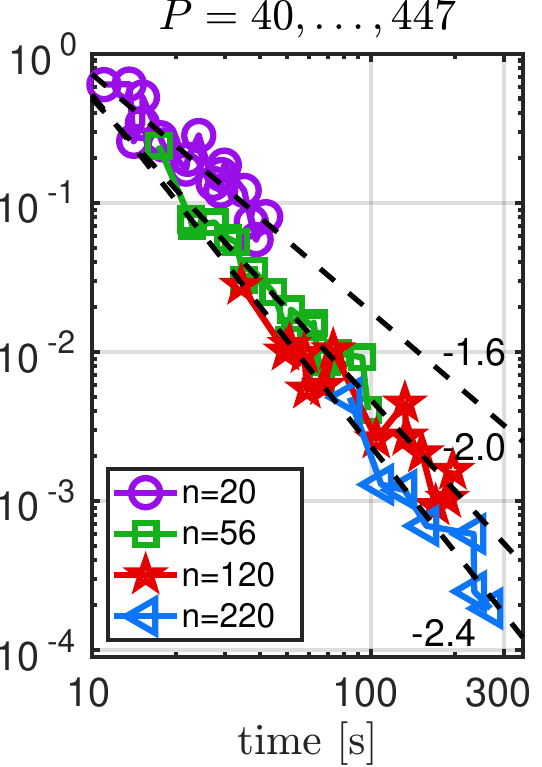}
  \caption{The relative maximum error for the 3D problem is shown for different numbers of local points $n$ as a function of the fill distance $h$ (left subplot) and as a function of the run time in seconds (right subplot) using $q=5$, $\varepsilon=1.2$, $H_s=1$, and $\delta_0=0.1$.}
  \label{fig:conv3D}
\end{figure}

Convergence and efficiency for the 3D problem are shown in \cref{fig:conv3D}. The theoretical convergence rates from~\eqref{eq:theorconv} here become 0.5, 2.5, 4.5, and 6.5. Note that here the amount of RAM in the laptop limits how large problems can be solved. Therefore the curves for larger $n$ have fewer points.
The values of $P_0$ that were used are 40, 60:10:90, 100:20:300, 360, 400.
The numerical convergence rates are also here better than expected, but in a less systematic way. We would need to run a larger range of problem sizes to really see the trends. The behavior with respect to efficiency is similar to that for the 2D case. The number of local points $n$ should increase when the error tolerance is lowered.

\begin{figure}[!htb]
  \centering
  \includegraphics[width=0.45\textwidth]{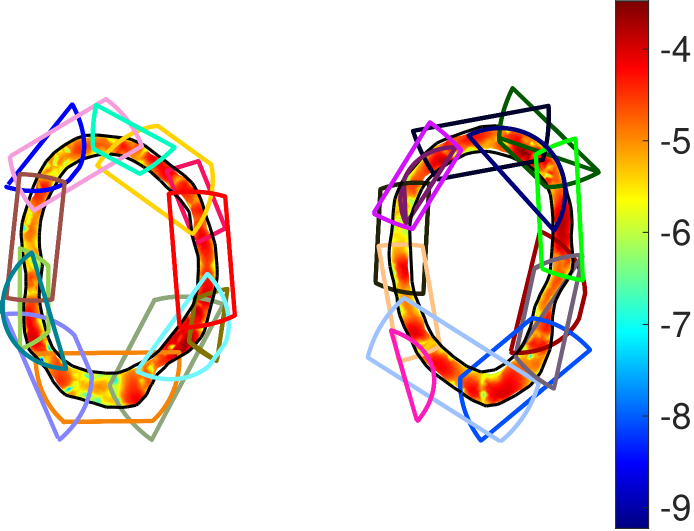}%
  \includegraphics[width=0.45\textwidth]{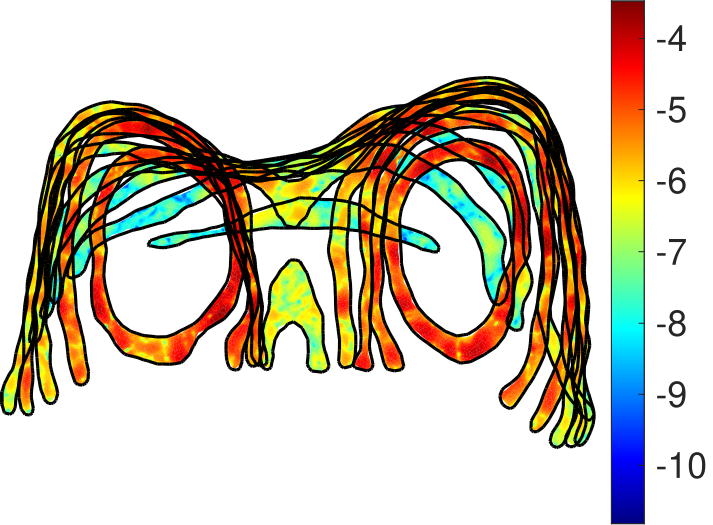}
  \caption{The logarithm in base 10 of the errors in one slice (left) and in ten slices (right) in the RBF-PUM solutions to the 3D problem for $n=220$ and $P=100$. In the left subfigure, also the intersections of the plane with the cover is shown.}\label{fig:slices}
\end{figure}

The errors in the 3D solution for one set of parameters is shown in different slices in \cref{fig:slices}. As can be seen in the left subplot, the patch structure in 3D becomes non-trivial, making it hard to say if and how the error is affected by the patch overlaps. However, the error distribution appears to be relatively even within each slice.

\section{Conclusions}\label{sec:conc}
In this paper, we have developed an geometry adaptive least squares RBF-PUM. We first used the method for reconstruction of a geometry from a noisy point cloud. The reconstruction is stable enough that it can produce good results for different parameter choices, and it does not produce spurious zero level sets. We showed how the parameters can be optimized using quality measures. This approach can be used for different types of thin geometries, but a similar method could also be defined for general surfaces, where the patch layer would then follow the surface, not the geometry.
In the second part of the paper we showed that we can achieve a sustained and predictable convergence behaviour in our reconstructed 2D and 3D geometries, in spite of the irregular patch and node layouts.


\bibliographystyle{siamplain}
\bibliography{refs}

\begin{thebibliography}{10}

\bibitem{ASLGG17}
{\sc S.~Aroudj, P.~Seemann, F.~Langguth, S.~Guthe, and M.~Goesele}, {\em
  Visibility-consistent thin surface reconstruction using multi-scale kernels},
  ACM Trans. Graph., 36 (2017), 187 (13~pages),
  \url{https://doi.org/10.1145/3130800.3130851}.

\bibitem{invive}
{\sc N.~Cacciani, E.~Larsson, A.~Lauro, M.~Meggiolaro, A.~Scatto, I.~Tominec,
  and P.-F. Villard}, {\em A first meshless approach to simulation of the
  elastic behaviour of the diaphragm}, in Spectral and High Order Methods for
  Partial Differential Equations ICOSAHOM 2018, S.~J. Sherwin, D.~Moxey,
  J.~Peir{\'o}, P.~E. Vincent, and C.~Schwab, eds., vol.~134 of Lect. Notes
  Comput. Sci. Eng., Springer International Publishing, Cham, 2020,
  pp.~501--512, \url{https://doi.org/10.1007/978-3-030-39647-3_40}.

\bibitem{Carr01}
{\sc J.~C. Carr, R.~K. Beatson, J.~B. Cherrie, T.~J. Mitchell, W.~R. Fright,
  B.~C. McCallum, and T.~R. Evans}, {\em Reconstruction and representation of
  3{D} objects with radial basis functions}, in SIGGRAPH '01: Proceedings of
  the 28th Annual Conference on Computer Graphics and Interactive Techniques,
  New York, NY, USA, 2001, Association for Computing Machinery, pp.~67--76,
  \url{https://doi.org/10.1145/383259.383266}.

\bibitem{Carr03}
{\sc J.~C. Carr, R.~K. Beatson, B.~C. McCallum, W.~R. Fright, T.~J. McLennan,
  and T.~J. Mitchell}, {\em Smooth surface reconstruction from noisy range
  data}, in Proceedings of the 1st International Conference on Computer
  Graphics and Interactive Techniques in Australasia and South East Asia,
  GRAPHITE '03, New York, NY, USA, 2003, Association for Computing Machinery,
  pp.~119--126, \url{https://doi.org/10.1145/604471.604495}.

\bibitem{CaLaMoMo06}
{\sc G.~Casciola, D.~Lazzaro, L.~Montefusco, and S.~Morigi}, {\em Shape
  preserving surface reconstruction using locally anisotropic radial basis
  function interpolants}, Comput. Math. Appl., 51 (2006), pp.~1185--1198,
  \url{https://doi.org/10.1016/j.camwa.2006.04.002}.

\bibitem{Cav20}
{\sc R.~Cavoretto and A.~De~Rossi}, {\em Error indicators and refinement
  strategies for solving {P}oisson problems through a {RBF} partition of unity
  collocation scheme}, Appl. Math. Comput., 369 (2020), pp.~124824, 18,
  \url{https://doi.org/10.1016/j.amc.2019.124824},
  \url{https://doi.org/10.1016/j.amc.2019.124824}.

\bibitem{Cav19}
{\sc R.~Cavoretto, A.~De~Rossi, G.~E. Fasshauer, M.~J. McCourt, and
  E.~Perracchione}, {\em Anisotropic weights for {RBF}-{PU} interpolation with
  subdomains of variable shapes}, in Numerical mathematics and advanced
  applications---{ENUMATH} 2017, F.~A. Radu, K.~Kumar, I.~Berre, J.~M.
  Nordbotten, and I.~S. Pop, eds., vol.~126 of Lect. Notes Comput. Sci. Eng.,
  Springer, Cham, 2019, pp.~93--101,
  \url{https://doi.org/10.1007/978-3-319-96415-7\_6}.

\bibitem{CheShch18}
{\sc G.~Cheng and V.~Shcherbakov}, {\em Anisotropic radial basis function
  methods for continental size ice sheet simulations}, J. Comput. Phys., 372
  (2018), pp.~161--177, \url{https://doi.org/10.1016/j.jcp.2018.06.020},
  \url{https://doi.org/10.1016/j.jcp.2018.06.020}.

\bibitem{DraFuWri21}
{\sc K.~P. Drake, E.~J. Fuselier, and G.~B. Wright}, {\em Implicit surface
  reconstruction with a curl-free radial basis function partition of unity
  method}, 2021, \url{https://arxiv.org/abs/2101.05940}.

\bibitem{Fass07}
{\sc G.~E. Fasshauer}, {\em Meshfree approximation methods with {MATLAB}},
  vol.~6 of Interdisciplinary Mathematical Sciences, World Scientific
  Publishing Co. Pte. Ltd., Hackensack, NJ, 2007,
  \url{https://doi.org/10.1142/6437}, \url{http://dx.doi.org/10.1142/6437}.

\bibitem{FassMc12}
{\sc G.~E. Fasshauer and M.~J. McCourt}, {\em Stable evaluation of {G}aussian
  radial basis function interpolants}, SIAM J. Sci. Comput., 34 (2012),
  pp.~A737--A762, \url{https://doi.org/10.1137/110824784},
  \url{https://doi.org/10.1137/110824784}.

\bibitem{FoFly15book}
{\sc B.~Fornberg and N.~Flyer}, {\em A primer on radial basis functions with
  applications to the geosciences}, vol.~87 of CBMS-NSF Regional Conference
  Series in Applied Mathematics, Society for Industrial and Applied Mathematics
  (SIAM), Philadelphia, PA, 2015,
  \url{https://doi.org/10.1137/1.9781611974041.ch1}.

\bibitem{FoLaFly11}
{\sc B.~Fornberg, E.~Larsson, and N.~Flyer}, {\em Stable computations with
  {G}aussian radial basis functions}, SIAM J. Sci. Comput., 33 (2011),
  pp.~869--892, \url{https://doi.org/10.1137/09076756X},
  \url{https://doi.org/10.1137/09076756X}.

\bibitem{FoPi07}
{\sc B.~Fornberg and C.~Piret}, {\em A stable algorithm for flat radial basis
  functions on a sphere}, SIAM J. Sci. Comput., 30 (2007), pp.~60--80,
  \url{https://doi.org/10.1137/060671991}.

\bibitem{FoWri04}
{\sc B.~Fornberg and G.~Wright}, {\em Stable computation of multiquadric
  interpolants for all values of the shape parameter}, Comput. Math. Appl., 48
  (2004), pp.~853--867, \url{https://doi.org/10.1016/j.camwa.2003.08.010}.

\bibitem{Halton60}
{\sc J.~H. Halton}, {\em On the efficiency of certain quasi-random sequences of
  points in evaluating multi-dimensional integrals}, Numer. Math., 2 (1960),
  pp.~84--90, \url{https://doi.org/10.1007/BF01386213},
  \url{https://doi.org/10.1007/BF01386213}.

\bibitem{invive_project}
{\sc INVIVE}, {\em The invive project}, 2017,
  \url{https://www.it.uu.se/research/scientific_computing/project/rbf/biomech}.

\bibitem{KoLaYu19}
{\sc K.~Kormann, C.~Lasser, and A.~Yurova}, {\em Stable interpolation with
  isotropic and anisotropic gaussians using hermite generating function}, SIAM
  J. Sci. Comput., 41 (2019), pp.~A3839--A3859,
  \url{https://doi.org/10.1137/19M1262449},
  \url{https://doi.org/10.1137/19M1262449}.

\bibitem{LaFo05}
{\sc E.~Larsson and B.~Fornberg}, {\em Theoretical and computational aspects of
  multivariate interpolation with increasingly flat radial basis functions},
  Comput. Math. Appl., 49 (2005), pp.~103--130,
  \url{https://doi.org/10.1016/j.camwa.2005.01.010}.

\bibitem{LLHF13}
{\sc E.~Larsson, E.~Lehto, A.~Heryudono, and B.~Fornberg}, {\em Stable
  computation of differentiation matrices and scattered node stencils based on
  {G}aussian radial basis functions}, SIAM J. Sci. Comput., 35 (2013),
  pp.~A2096--A2119, \url{https://doi.org/10.1137/120899108},
  \url{https://doi.org/10.1137/120899108}.

\bibitem{LaMaMiPoo22}
{\sc E.~Larsson, B.~Mavri\v{c}, A.~Michael, and F.~Pooladi}, {\em A numerical
  investigation of some {RBF}-{FD} error estimates}, Dolomites Res. Notes
  Approx., 15 (2022), pp.~78--95.

\bibitem{RBFQR3Dmixed}
{\sc E.~Larsson and A.~Michael}, {\em Computing mixed derivatives in {3D} with
  the {RBF}-{QR} method}.
\newblock manuscript in preparation, 2024.

\bibitem{LaShchHer17}
{\sc E.~Larsson, V.~Shcherbakov, and A.~Heryudono}, {\em A least squares radial
  basis function partition of unity method for solving {PDE}s}, SIAM J. Sci.
  Comput., 39 (2017), pp.~A2538--A2563,
  \url{https://doi.org/10.1137/17M1118087}.

\bibitem{levy2010}
{\sc B.~L{\'e}vy and Y.~Liu}, {\em ${L}_p$ centroidal {V}oronoi tessellation
  and its applications}, ACM Trans. Graphics, 29 (2010), 119 (11~pages),
  \url{https://doi.org/10.1145/1778765.1778856}.

\bibitem{LiWaBruWa16}
{\sc S.~Liu, C.~C. Wang, G.~Brunnett, and J.~Wang}, {\em A closed-form
  formulation of {HRBF}-based surface reconstruction by approximate solution},
  Comput.-Aided Des., 78 (2016), pp.~147--157,
  \url{https://doi.org/10.1016/j.cad.2016.05.001}.

\bibitem{Micchelli86}
{\sc C.~A. Micchelli}, {\em Interpolation of scattered data: distance matrices
  and conditionally positive definite functions}, Constr. Approx., 2 (1986),
  pp.~11--22, \url{https://doi.org/10.1007/BF01893414}.

\bibitem{OhBeSe06}
{\sc Y.~Ohtake, A.~Belyaev, and H.-P. Seidel}, {\em Sparse surface
  reconstruction with adaptive partition of unity and radial basis functions},
  Graph. Models, 68 (2006), pp.~15--24,
  \url{https://doi.org/10.1016/j.gmod.2005.08.001}.

\bibitem{PaSka11}
{\sc R.~Pan and V.~Skala}, {\em A two-level approach to implicit surface
  modeling with compactly supported radial basis functions}, Eng. Comput., 27
  (2011), pp.~299--307, \url{https://doi.org/10.1007/s00366-010-0199-1}.

\bibitem{Piret12}
{\sc C.~Piret}, {\em The orthogonal gradients method: a radial basis functions
  method for solving partial differential equations on arbitrary surfaces}, J.
  Comput. Phys., 231 (2012), pp.~4662--4675,
  \url{https://doi.org/10.1016/j.jcp.2012.03.007}.

\bibitem{PiDu16}
{\sc C.~Piret and J.~Dunn}, {\em Fast {RBF} {OG}r for solving {PDE}s on
  arbitrary surfaces}, in Numerical Computations: Theory and Algorithms
  (NUMTA-2016), Y.~D. Sergeyev, D.~E. Kvasov, F.~Dell'Accio, and M.~S.
  Mukhametzhanov, eds., vol.~1776 of AIP Conf. Proc., Amer. Inst. Phys.,
  Melville, NY, 2016.

\bibitem{RieZwi10}
{\sc C.~Rieger and B.~Zwicknagl}, {\em Sampling inequalities for infinitely
  smooth functions, with applications to interpolation and machine learning},
  Adv. Comput. Math., 32 (2010), pp.~103--129,
  \url{https://doi.org/10.1007/s10444-008-9089-0}.

\bibitem{SVHL15}
{\sc A.~Safdari-Vaighani, A.~Heryudono, and E.~Larsson}, {\em A radial basis
  function partition of unity collocation method for convection-diffusion
  equations arising in financial applications}, J. Sci. Comput., 64 (2015),
  pp.~341--367, \url{https://doi.org/10.1007/s10915-014-9935-9},
  \url{https://doi.org/10.1007/s10915-014-9935-9}.

\bibitem{Schaback05}
{\sc R.~Schaback}, {\em Multivariate interpolation by polynomials and radial
  basis functions}, Constr. Approx., 21 (2005), pp.~293--317,
  \url{https://doi.org/10.1007/s00365-004-0585-2}.

\bibitem{Scho38}
{\sc I.~J. Schoenberg}, {\em Metric spaces and completely monotone functions},
  Ann. of Math. (2), 39 (1938), pp.~811--841,
  \url{https://doi.org/10.2307/1968466}.

\bibitem{Shepard68}
{\sc D.~Shepard}, {\em A two-dimensional interpolation function for
  irregularly-spaced data}, in Proceedings of the 1968 23rd ACM National
  Conference, ACM '68, New York, NY, USA, 1968, ACM, pp.~517--524,
  \url{https://doi.org/10.1145/800186.810616}.

\bibitem{ToBre21}
{\sc I.~Tominec and E.~Breznik}, {\em An unfitted {RBF}-{FD} method in a
  least-squares setting for elliptic {PDE}s on complex geometries}, J. Comput.
  Phys., 436 (2021), 110283 (24~pages),
  \url{https://doi.org/10.1016/j.jcp.2021.110283}.

\bibitem{ToLaHe21}
{\sc I.~Tominec, E.~Larsson, and A.~Heryudono}, {\em A least squares radial
  basis function finite difference method with improved stability properties},
  SIAM J. Sci. Comput., 43 (2021), pp.~A1441--A1471,
  \url{https://doi.org/10.1137/20M1320079}.

\bibitem{TVLBC21}
{\sc I.~Tominec, P.~Villard, E.~Larsson, V.~Bayona, and N.~Cacciani}, {\em An
  unfitted radial basis function generated finite difference method applied to
  thoracic diaphragm simulations}, 2021,
  \url{https://arxiv.org/abs/2103.03673},
  \url{https://arxiv.org/abs/2103.03673}.

\bibitem{villard11}
{\sc P.-F. Villard, P.~Boshier, F.~Bello, and D.~Gould}, {\em Virtual reality
  simulation of liver biopsy with a respiratory component}, in Liver Biopsy,
  H.~Takahashi, ed., InTechOpen, 2011, \url{https://doi.org/10.5772/22033}.

\bibitem{Wend95}
{\sc H.~Wendland}, {\em Piecewise polynomial, positive definite and compactly
  supported radial functions of minimal degree}, Adv. Comput. Math., 4 (1995),
  pp.~389--396, \url{https://doi.org/10.1007/BF02123482}.

\bibitem{Wend02}
{\sc H.~Wendland}, {\em Fast evaluation of radial basis functions: methods
  based on partition of unity}, in Approximation theory, {X} ({S}t. {L}ouis,
  {MO}, 2001), Innov. Appl. Math., Vanderbilt Univ. Press, Nashville, TN, 2002,
  pp.~473--483.

\bibitem{WMTPM21}
{\sc R.~M. Whebell, T.~J. Moroney, I.~W. Turner, R.~Pethiyagoda, and S.~W.
  McCue}, {\em Implicit reconstructions of thin leaf surfaces from large, noisy
  point clouds}, Appl. Math. Model., 98 (2021), pp.~416--434,
  \url{https://doi.org/10.1016/j.apm.2021.05.014}.

\bibitem{FoWri17}
{\sc G.~B. Wright and B.~Fornberg}, {\em Stable computations with flat radial
  basis functions using vector-valued rational approximations}, J. Comput.
  Phys., 331 (2017), pp.~137--156,
  \url{https://doi.org/10.1016/j.jcp.2016.11.030},
  \url{https://doi.org/10.1016/j.jcp.2016.11.030}.

\bibitem{QiWaWu06}
{\sc Q.~Xia, M.~Y. Wang, and X.~Wu}, {\em Orthogonal least squares in partition
  of unity surface reconstruction with radial basis function}, in Geometric
  Modeling and Imaging--New Trends (GMAI'06), 2006, pp.~28--33,
  \url{https://doi.org/10.1109/GMAI.2006.40}.

\bibitem{YWZP08}
{\sc J.~Yang, Z.~Wang, C.~Zhu, and Q.~Peng}, {\em Implicit surface
  reconstruction with radial basis functions}, in Computer Vision and Computer
  Graphics. Theory and Applications, J.~Braz, A.~Ranchordas, H.~J. Ara{\'u}jo,
  and J.~M. Pereira, eds., Springer Berlin Heidelberg, 2008, pp.~5--12,
  \url{https://doi.org/10.1007/978-3-540-89682-1_1}.

\end{thebibliography}
\end{document}